%
\magnification=\magstep1
\input amstex
\UseAMSsymbols
\input pictex
\NoBlackBoxes
  \font\gross=cmbx10 scaled\magstep1 
   \font\rmk=cmr8      \font\ttk=cmtt8

   \newcount\notenumber
   
   \def\note{\advance\notenumber by 1 
       \plainfootnote{$^{\the\notenumber}$}}  
\def\fin{\operatorname{fin}}
\def\Odim{\operatorname{Odim}}
\def\m{\bold m}

\def\mod{\operatorname{mod}}
\def\Max{\operatorname{Max}}

\def\Im{\operatorname{Im}}

\def\rep{\operatorname{rep}}

\def\Hom{\operatorname{Hom}}

\def\End{\operatorname{End}}
\def\Ext{\operatorname{Ext}}
\def\rad{\operatorname{rad}}
\def\soc{\operatorname{soc}}
\def\top{\operatorname{top}}
\def\add{\operatorname{add}}

\def\op{{\text{op}}}
\def\repdim{\operatorname{repdim}}
\def\LL{\operatorname{LL}}

\def\arr#1#2{\arrow <1.5mm> [0.25,0.75] from #1 to #2}

 
\def\vect#1#2#3#4{\beginpicture
   \setcoordinatesystem units <.15cm,.15cm>
   \put{$\ssize #1$} at 1 2
   \put{$\ssize #2$} at 0 1
   \put{$\ssize #3$} at 2 1
   \put{$\ssize #4$} at 1 0
   \endpicture}

\vglue2truecm

\centerline{\gross On the representation dimension of artin algebras.}
	\bigskip
\centerline{\bf Claus Michael Ringel}
	\bigskip\bigskip
{\narrower
{\bf Abstract.} The representation dimension of an artin algebra 
as introduced  by M.~Auslander
in his Queen Mary Notes is the minimal possible global dimension of
the endomorphism ring of a generator-cogenerator. The following report is based 
on two texts written in 2008 in connection with a workshop at Bielefeld. 
The first part presents a full proof that  
any torsionless-finite artin algebra has representation dimension at most $3$,
and provides a long list of classes of algebras which are torsionless-finite.
In the second part we show that the representation 
dimension is adjusted very well to forming tensor products of algebras.
In this way one obtains a wealth
of examples of artin algebras with large representation dimension.
In particular, we show: The tensor product of $n$ 
representation-infinite path algebras of 
bipartite quivers has representation dimension precisely $n+2$.\par}
	\bigskip\bigskip
Let $\Lambda$ be an artin algebra.
The representation dimension $\repdim \Lambda$ 
of $\Lambda$ was introduced 1971 by M.~Auslander
in his Queen Mary Notes [A], it is the minimal possible global dimension of
the endomorphism ring of a generator-cogenerator (a {\it generator-cogenerator}
is a $\Lambda$-module $M$ such that any indecomposable projective or injective
$\Lambda$-module is a direct summand of $M$); a generator-cogenerator $M$ such that
the global dimension of the endomorphism ring of $M$ is minimal will be said to
be an {\it Auslander generator.} All the classes of algebras where
Auslander was able to determine the precise representation dimension
turned out to have representation dimension at most $3$. Thus, he 
asked, on the one hand, whether the
representation dimension can be greater than $3$, but also, on the other hand, 
whether it always has to be finite. 
These questions have been answered only recently: The finiteness
of the representation dimension was shown by Iyama [I1] in 2003
(for a short proof using the notion of strongly quasi-hereditary algebras see [R4]). 
For some artin algebras $\Lambda$, one knows that $\repdim \Lambda \le \LL(\Lambda)+1,$ 
where $\LL(\Lambda)$ is the Loewy length of $\Lambda$.
On the other hand, Rouquier [Rq] has shown that the representation dimension of the exterior
algebra of a vector space of dimension $n \ge 1$ is $n+1.$ 

The name {\it representation dimension} was coined by Auslander on the basis of his
observation that $\Lambda$ is representation-finite if and only if $\repdim \Lambda \le 2$.
Not only Auslander, but later also many other mathematicians were able to prove
that the representation dimension of many well-known classes of artin algebras is
bounded by $3$. 
For artin algebras with representation dimension at most $3$, Igusa and Todorov
[IT] have shown that the finitistic dimension of these algebras is finite. 
	\medskip

Our report is based on two texts written in 2008 in connection with a
Bielefeld workshop on the representation dimension [Bi].

Part I is a modified version of [R3] which was written
as an introduction for the workshop, it aimed at a general scheme
for some of the known proofs for the upper bound $3$, by providing the assertion 4.2
that 
{\it any torsionless-finite artin algebra has representation dimension at most $3$}.
We recall that a $\Lambda$-module is said to be {\it torsionless} or {\it divisible}
provided it is a submodule of a projective module, or a factor module of an
injective module, respectively. The artin algebra $\Lambda$ is said to
be {\it torsionless-finite} provided there are only finitely many
isomorphism classes of indecomposable torsionless $\Lambda$-modules.
Two ingredients are needed for the proof of 4.2, one is
the characterization 4.3 of Auslander generators which was used already
by Auslander (at least implicitly)
in the Queen Mary notes, the second is the bijection 3.2 between the
isomorphism classes of the indecomposable torsionless and the indecomposable
divisible modules which can be found in the appendix of Auslander-Bridger [AB].

In Part II we want to outline that the representation 
dimension is adjusted very well to forming tensor products of algebras.
It should be stressed that already in 2000, Changchang Xi [X1] was
drawing the attention to this relationship by showing that
$$
 \repdim \Lambda\otimes_k \Lambda' \le \repdim \Lambda + \repdim \Lambda'
$$
for finite-dimensional $k$-algebras $\Lambda, \Lambda'$, provided 
$k$ is a perfect field. In 2009,
Oppermann  [O1] gave a lattice criterion for obtaining a lower
bound for the representation dimension (see 6.3) and we show (see 7.1)
that this lattice construction is compatible with tensor products; this
fact was also noted by Oppermann [O4]. In this way one obtains a wealth
of examples with large representation dimension. These sections 6 and 7 are 
after-thoughts to the workshop-lecture of Oppermann and a corresponding 
text was privately distributed after the workshop. 
The last two sections show that sometimes
one is able to determine the precise value of the representation dimension.
Namely, we will show in 9.4:
{\it Let $\Lambda_1,\dots\Lambda_n$ be representation-infinite path algebras of 
bipartite quivers. Then the algebra
$\Lambda = \Lambda_1\otimes_k\cdots\otimes_k\Lambda_n$
has representation dimension precisely $n+2$.} For quivers without multiple
arrows, this result was presented at the Abel conference in Balestrand, June 2011,
the example of the tensor product of two copies of the Kronecker algebra was exhibited
already by Oppermann at the Bielefeld workshop. 
	\bigskip\bigskip
We consider an artin algebra $\Lambda$ with duality functor $D$. Usually,
we will consider left $\Lambda$-modules of finite length and call them just
modules. Always, morphisms will be written on the opposite side of the scalars.
	\medskip
Given a class $\Cal M$ of modules, we denote by $\add \Cal M$ the modules which are
(isomorphic to) direct summands of direct sums of modules in $\Cal M$. If $M$
is a module, we write $\add M = \add \{M\}.$
We say that $\Cal M$
is {\it finite} provided there are only finitely many isomorphism classes of indecomposable
modules in $\add \Cal M$, thus provided there exists a module $M$ with 
$\add \Cal M = \add M.$
	\bigskip
{\bf Acknowledgment.} The author is indebted to Dieter Happel and Steffen Oppermann
for remarks concerning the presentation of the paper. 

	\vfill\eject
\centerline{\gross Part I. Torsionless-finite artin algebras.}
	\bigskip
{\bf 1. The torsionless modules for $\Lambda$ and $\Lambda^\op$} 
	\medskip
Let $\Cal L= \Cal L(\Lambda)$ be the class of torsionless $\Lambda$-modules
and $\Cal P= \Cal P(\Lambda)$ the class of projective $\Lambda$-modules. 
Then $\Cal P(\Lambda) \subseteq \Cal L(\Lambda),$ and we may consider
 the factor category 
$\Cal L(\Lambda)/\Cal P(\Lambda)$ obtained from
$\Cal L(\Lambda)$ by factoring out the ideal of all maps which factor
through a projective module. 
	\bigskip
{\bf (1.1) Theorem.} {\it There is a duality 
$$
 \eta\: \Cal L(\Lambda)/\Cal P(\Lambda) \longrightarrow 
 \Cal L(\Lambda^\op)/\Cal P(\Lambda^\op)
$$
with the following property: If $U$ is a torsionless module, and
$f\:P_1(U) \to P_0(U)$ is a projective presentation of $U$, then
for $\eta(U)$ we can take the image of $\Hom(f,\Lambda)$.}
	\medskip
In the proof, we will use the following definition: 
We call an exact sequence $P_1 \to P_0 \to P_{-1}$
with projective modules $P_i$ {\it strongly exact} provided it remains exact
when we apply $\Hom(-,\Lambda).$ Let $\Cal E$ be the category of strongly exact
sequences $P_1 \to P_0 \to P_{-1}$
with projective modules $P_i$ (as a full subcategory of the category of complexes).
	\medskip
{\bf (1.2) Lemma.} {\it The exact sequence $ P_1 @>f>> P_0 @>g>> P_{-1},$ with 
all $P_i$ projective and epi-mono factorization $g = ue$ 
is strongly exact if and only if $u$ is a left $\Lambda$-approximation.}
	\medskip
Proof: Under the functor $\Hom(-,\Lambda)$, we obtain
$$
 \Hom(P_{-1},\Lambda) @>g^*>> \Hom(P_{0},\Lambda) @>f^*>> \Hom(P_{1},\Lambda)
$$
with zero composition. Assume that $u$ is a left $\Lambda$-approximation.
Given $\alpha\in \Hom(P_{0},\Lambda)$ with $f^*(\alpha) = 0,$ we rewrite
$f^*(\alpha) = \alpha f.$ Now $e$ is a cokernel of $f$, thus there is $\alpha'$
with $\alpha = \alpha'e.$ Since $u$ is a left $\Lambda$-approximation, there is
$\alpha''$ with $\alpha' = \alpha''u.$ It follows that $\alpha = \alpha'e = \alpha''ue
= \alpha''g  = g^*(\alpha'').$ 

Conversely, assume that the sequence $(*)$ is exact, let $U$ be the image of $g$, thus
$e\:P_0 \to U, u\:U \to P_{-1}.$ Consider a map $\beta\:U \to \Lambda$. Then 
$f^*(\beta e) = \beta e f = 0,$ thus there is $\beta'\in \Hom(P_{-1},\Lambda$
with $g^*(\beta') = \beta e.$ But $g^*(\beta) = \beta' g = \beta' ue$ and $\beta e =
\beta' ue$ implies $\beta = \beta'u$, since $e$ is an epimorphism.
	\medskip
{\bf Proof of Theorem 1.1.}
Let $\Cal U$ be the full subcategory of $\Cal E$ of all sequences which are direct sums
of sequences of the form
$$
 P @>>> 0 @>>> 0, \quad 
 P @>1>> P @>>> 0, \quad 
 0 @>>> P @>1>> P, \quad 
 0 @>>> 0 @>>> P.
$$
In order to define the functor $q\:\Cal E \to \Cal L$, let $q(P_1 @>f>> P_0 @>g>> P_{-1})$
be the image of $g.$
Clearly, $q$ sends $\Cal U$ onto $\Cal P$, thus it induces a functor
$$
 \overline q\: \Cal E/\Cal U \longrightarrow \Cal L/\Cal P.
$$
Claim: {\it This functor $\overline q$ is an equivalence.}
	\medskip
First of all, the functor $q$ is dense: starting with $U\in \Cal L,$ let 
$$
 P_1 @>f>> P_0 @>e>> U @>>> 0
$$
be a projective presentation of $U$, let $u\:U \to P_{-1}$ be a left 
$\Lambda$-approximation of $U$, and $g = ue.$

Second, the functor $q$ is full. This follows from the lifting properties of
projective presentations and left $\Lambda$-approximations.

It remains to show that $\overline q$ is faithful. We will give the proof in detail
(and it may look quite technical), however we should remark that all the
arguments are standard; they are the usual ones dealing with homotopy categories of
complexes. Looking at strongly exact sequences $ P_1 @>f>> P_0 @>g>> P_{-1}$,
one should observe that the image $U$ of $g$ has to be considered as the essential
information: starting from $U$, one may attach to it a projective presentation
(this means going from $U$ to the left in order to obtain $P_1 @>f>> P_0$)
as well as a left $\Lambda$-approximation of $U$  
(this means going from $U$ to the right in order to obtain $P_{-1}$).
	\medskip
In order to show that $\overline q$ is faithful, let us consider 
the following commutative diagram
$$
\CD
  P_1   @>f>>   P_0   @>g>>   P_{-1} \cr
 @Vh_1VV     @Vh_0VV       @Vh_{-1}VV  \cr
  P'_1  @>f'>>   P'_0  @>g'>>   P'_{-1}
\endCD
$$
with strongly exact rows. We consider epi-mono factorizations $g = eu, g' = e'u'$ with
$e\:P_0 \to U, u\:U \to P_{-1}, e'\:P'_0 \to U', u'\:U' \to P'_{-1},$ thus
$q(P_\bullet) = U, q(P'_\bullet) = U'$. Assume that $q(h_\bullet) = ab,$ where
$a\:U \to X, b\:X \to U'$ with $X$ projective. We have to show that $h_\bullet$
belongs to $\Cal U.$

The factorizations $g = eu, g' = e'u', q(h_\bullet) = ab$ provide the
following equalities: 
$$
   eab = h_0e',\quad uh_1 = abu'.
$$
Since $u\:U \to P_{-1}$ is a left $\Lambda$-approximation and $X$ is projective,
there is $a'\:P_{-1} \to X$ with $ua' = a$. Since $e'\:P'_0\to U'$ 
is surjective and $X$ is projective, there is $b'\:X \to P'_0$ with $b'e' = b$.

Finally, we need $c\:P_0 \to P'_1$ with $cf' = h_0-eab'$. 
Write $f' = w'v'$ with $w'$ epi and $v'$ mono; in particular, $v'$ is the kernel
of $g'$. Note that
$eab'g' = eab'e'u'= eabu' = h_0e'u' = h_0g'$, thus 
$(h_0-eab')g' = h_0g'-eab'g' = h_0g'-h_0g' = 0,$ thus $h_0-eab'$ factors through
the kernel $v'$ of $g'$, say $h_0-eab' = c'v'$. Since $P_0$ is projective and $w'$
is surjective, we find $c\:P_0 \to P'_1$ with $cw' = c'$, thus $cf' = cw'v' =
c'v' = h_0-eab'.$ 

Altogether, we obtain the following commutative diagram
$$
\CD
  P_1   @>f>>   P_0   @>g>>   P_{-1} \cr
 @V{\left[\smallmatrix 1 & f\endsmallmatrix\right]}VV     
   @VV{\left[\smallmatrix 1 & ea \endsmallmatrix\right]}V       
   @VV{\left[\smallmatrix a' & h_1\!-\!a'bu'\endsmallmatrix\right]}V  \cr
  P_1\oplus P_0   @>{\left[\smallmatrix 0 & 0 \cr 1 & 0\endsmallmatrix\right]}>>   
                P_0\oplus X   
                  @>{\left[\smallmatrix 0 & 0 \cr 1 & 0\endsmallmatrix\right]}>>                   X\oplus P'_{-1} \cr
 @V{\left[\smallmatrix h_1-fc\strut \cr\strut c \endsmallmatrix\right]}VV     
    @VV{\left[\smallmatrix h_0-eab'\strut \cr \strut b' \endsmallmatrix\right]}V
    @VV{\left[\smallmatrix bu' \strut\cr\strut 1\endsmallmatrix\right]} V  \cr
  P'_1  @>f'>>   P'_0  @>g'>>   P'_{-1}
\endCD
$$
which is the required factorization of $h_\bullet$ (here, the commutativity
of the four square has to be checked; in addition, one has to verify
that the vertical compositions yield the maps $h_i$; all these calculations
are straight forward). 
	\medskip
Now consider the functor $\Hom(-,\Lambda)$, it yields a duality
$$
 \Hom(-,\Lambda)\: \Cal E(\Lambda) \longrightarrow \Cal E(\Lambda^\op)
$$
which sends $\Cal U(\Lambda)$ onto $\Cal U(\Lambda^\op).$ Thus, we obtain a 
duality
$$
  \Cal E(\Lambda)/\Cal U(\Lambda) \longrightarrow 
  \Cal E(\Lambda^\op)/\Cal U(\Lambda^\op).
$$

Combining the functors considered, we obtain the following sequence
$$
 \Cal L(\Lambda)/\Cal P(\Lambda) @<\;\;\overline q\;\;<< \Cal E(\Lambda)/\Cal U(\Lambda)
 @>\Hom(-,\Lambda)>> \Cal E(\Lambda^\op)/\Cal U(\Lambda^\op) @>\;\;\overline q\;\;>>  \Cal L(\Lambda^\op)/\Cal P(\Lambda^\op),
$$
this is duality, and we denote it by $\eta$.

It remains to show that $\eta$ is given by the mentioned recipe.
Thus, let $U$ be a torsionless module. Take a projective presentation 
$$
 P_1 @>f>> P_0 @>e>> U @>>> 0
$$
of $U$, and let $m\:U \to P_{-1}$ be a left $\Cal P$-approximation of $U$ and $g = eu.$
Then 
$$
 P_\bullet = (P_1 @>f>> P_0 @>g>> P_{-1})
$$ 
belongs to $\Cal E$ and 
$q(P_\bullet) = U.$ The functor $\Hom(-,\Lambda)$ sends $P_\bullet$ to
$$
 \Hom(P_\bullet,\Lambda) = 
 (\Hom(P_{-1},\Lambda) @>\Hom(g,\Lambda)>> \Hom(P_0,\Lambda) @>\Hom(f,\Lambda)>> \Hom(P_1,\Lambda))
$$
in $\Cal E(\Lambda^\op)$. Finally, the equivalence
$$
 \Cal E(\Lambda^\op)/\Cal U(\Lambda^\op) @>\;\;\overline q\;\;>>  \Cal L(\Lambda^\op)/\Cal P(\Lambda^\op)
$$
sends $\Hom(P_\bullet,\Lambda)$ to the image of $\Hom(f,\Lambda)$. 
	\bigskip\bigskip
{\bf 2. Consequences} 
	\medskip
{\bf (2.1) Corollary.} {\it There is a canonical bijection between the isomorphism
classes of the indecomposable torsionless $\Lambda$-modules and the 
isomorphism classes of the 
indecomposable torsionless $\Lambda^\op$-modules.}
	\medskip
Proof: The functor $\Hom(-,\Lambda)$ provides a bijection between the 
isomorphism
classes of the indecomposable projective $\Lambda$-modules and the 
isomorphism classes of the 
indecomposable projective $\Lambda^\op$-modules. For the non-projective
indecomposable torsionless modules, we use the duality $\eta$ given by
Theorem 1.
	\medskip
Remark. As we have seen, there are canonical bijections between
the indecomposable projective $\Lambda$-modules and 
and the indecomposable projective $\Lambda^\op$-modules,
as well between 
the indecomposable non-projective torsionless 
$\Lambda$-modules and the indecomposable non-projective torsionless 
$\Lambda^\op$-modules, both bijections being given by categorical
dualities, but these bijections do not combine to a bijection 
with nice categorical properties. We will exhibit suitable examples below. 
	\bigskip
{\bf (2.2) Corollary.} {\it If $\Lambda$ is torsionless-finite, also $\Lambda^\op$
is torsionless-finite.}
	 \bigskip
Whereas corollaries 2.1 and 2.2 are of interest only for 
non-commutative artin algebras,
the theorem itself is also of interest for $\Lambda$ commutative. 
	\medskip
{\bf (2.3) Corollary.} {\it For $\Lambda$ a commutative artin algebra, the category 
$\Cal L/\Cal P$ has a  self-duality.}
	\medskip
For example, consider the factor algebra 
$\Lambda = k[T]/\langle T^n\rangle$ of the polynomial ring $k[T]$ in one
variable, with $k$ is a field. Since $\Lambda$ is self-injective, all the modules
are torsionless. Note that in this case, $\eta$ coincides with the syzygy functor
$\Omega.$ 
	\bigskip\bigskip
{\bf 3. The torsionless and the divisible $\Lambda$-modules} 
	\medskip
Let $\Cal K = \Cal K(\Lambda)$ be the class of divisible $\Lambda$-modules. Of course, the duality
functor $D$ provides a bijection between the isomorphism classes of
divisible modules and the isomorphism classes of torsionless right modules.

We denote by $\Cal Q= \Cal Q(\Lambda)$ 
the class of injective modules. Clearly,  $D$ provides a duality
$$
 D\:  \Cal L(\Lambda^\op)/\Cal P(\Lambda^\op) \longrightarrow 
 \Cal K(\Lambda)/\Cal Q(\Lambda).
$$
Thus, we can reformulate theorem 1 as follows: 
{\it The categories $\Cal L(\Lambda)/\Cal P(\Lambda)$
and $\Cal K(\Lambda)/\Cal Q(\Lambda)$ are equivalent under the functor $D\eta$.} It
seems to be worthwhile to replace the functor $D\eta$ by the functor $\Sigma\tau$.
Here, $\tau$ is the Auslander-Reiten translation and $\Sigma$ is the suspension
functor (defined by 
$\Sigma(V) = I(V)/V,$ where $I(V)$ is an injective envelope of $V$).
Namely, in order to calculate $\tau(U)$, we start with a minimal projective
presentation $f\:P_1 \to P_0$ and take as $\tau(U)$ the kernel of 
$$
 D\Hom(f,\Lambda)\:D\Hom(P_1,\Lambda) \longrightarrow D\Hom(P_0,\Lambda).
$$
Now the kernel inclusion $\tau(U) \subset D\Hom(P_1,\Lambda)$ is an injective envelope
of $\tau(U);$ thus $\Sigma\tau(U)$ is the image of $D\Hom(f,\Lambda)$, but this
image is $D\eta(U).$ Thus we see that Theorem 1.1
can be formulated also as follows: 
	\bigskip
{\bf (3.1) Theorem.} {\it The categories $\Cal L(\Lambda)/\Cal P(\Lambda)$
and $\Cal K(\Lambda)/\Cal Q(\Lambda)$ are equivalent under the functor 
$\gamma = \Sigma\tau$.}
	\bigskip

{\bf (3.2) Corollary.} {\it If $\Lambda$ is torsionless-finite, the number
of isomorphism classes of indecomposable divisible modules is
equal to the number of isomorphism classes of indecomposable torsionless modules.}
	\bigskip
{\bf (3.3) Examples.} We insert here four examples so that one may 
get a feeling about the bijection 
between the isomorphism classes of indecomposable torsionless modules and those
of the indecomposable divisible modules.
	\smallskip
(1) The path algebra of a linearly oriented quiver of type
$A_3$ modulo the square of its radical.
$$
{\beginpicture
\setcoordinatesystem units <1cm,1cm>
\multiput{$\circ$} at 0 0  1 0  2 0 /
\arr{0.8 0}{0.2 0}
\arr{1.8 0}{1.2 0}
\setdots <.9mm>
\setquadratic
\plot 0.5 0.1  1 0.3  1.5 0.1 /
\endpicture}
$$
We present twice the Auslander-Reiten quiver. Left, we mark by $+$ the
indecomposable torsionless modules and encircle the unique non-projective torsionless
module. On the right, we mark by $*$ the indecomposable divisible
modules and encircle the unique non-injective divisible module:
$$
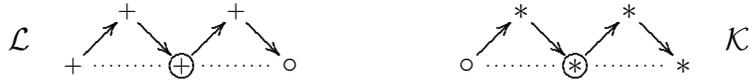
{\beginpicture
\setcoordinatesystem units <1cm,1cm>
\put{\beginpicture
\setcoordinatesystem units <.6cm,.6cm>
\put{$\circ$} at 4 0
\multiput{$\ssize+$} at 0 0  1 1  2 0  3 1  /
\circulararc 360 degrees from 2.22 0 center at 2 0 
\arr{0.2 0.2}{0.8 0.8}
\arr{1.2 0.8}{1.8 0.2}
\arr{2.2 0.2}{2.8 0.8}
\arr{3.2 0.8}{3.8 0.2}
\setdots <.9mm>
\plot 0.4 0  1.6 0 /
\plot 2.4 0  3.6 0 /
\put{$\Cal L$} at -1 0.5 
\endpicture} at 0 0
\put{\beginpicture
\setcoordinatesystem units <.6cm,.6cm>
\put{$\circ$} at 0 0
\multiput{$*$} at 1 1  2 0  3 1  4 0  /
\circulararc 360 degrees from 2.22 0 center at 2 0 
\arr{0.2 0.2}{0.8 0.8}
\arr{1.2 0.8}{1.8 0.2}
\arr{2.2 0.2}{2.8 0.8}
\arr{3.2 0.8}{3.8 0.2}
\setdots <.9mm>
\plot 0.4 0  1.6 0 /
\plot 2.4 0  3.6 0 /
\put{$\Cal K$} at 5 0.5 
\endpicture} at 5 0
\endpicture}
$$
	\medskip
(2) Next, we look at the algebra $\Lambda$ given by 
the following quiver with a commutative square; to the right, we present
its Auslander-Reiten quiver $\Gamma(\Lambda)$ 
and mark the torsionless and divisible modules
as in the previous example. Note that the subcategories $\Cal L$ and $\Cal K$
are linearizations of posets.
$$
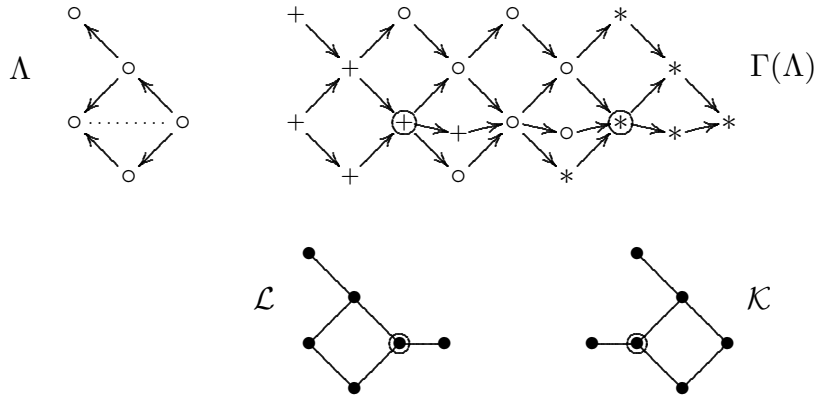
{\beginpicture
\setcoordinatesystem units <1cm,1cm>

\put{\beginpicture
\setcoordinatesystem units <.6cm,.6cm>
\multiput{$\circ$} at 0 1  0 3  1 0  1 2  2 1 /
\arr{0.8 1.8}{0.2 1.2}
\arr{0.8 2.2}{0.2 2.8}
\arr{0.8 0.2}{0.2 0.8}
\arr{1.8 0.8}{1.2 0.2}
\arr{1.8 1.2}{1.2 1.8}
\setdots <1mm>
\plot 0.3 1  1.7 1 /
\put{$\Lambda$} at -1 2
\endpicture} at 0 0
\put{\beginpicture
\setcoordinatesystem units <.6cm,.6cm>
\multiput{$\circ$} at 2 3  3 0  4 1  4 3  5 0.8  3 2  5 2 /
\multiput{$\ssize+$} at 0 1  0 3  1 0  1 2  2 1  3 0.8  /
\multiput{$*$} at 5 0    6 1  6 3  7 0.8  7 2  8 1  /
\circulararc 360 degrees from 2.22 1 center at 2 1 
\circulararc 360 degrees from 6.22 1 center at 6 1 
\arr{0.2 1.2}{0.8 1.8}
\arr{2.2 1.2}{2.8 1.8}
\arr{4.2 1.2}{4.8 1.8}
\arr{6.2 1.2}{6.8 1.8}

\arr{1.2 0.2}{1.8 0.8}
\arr{3.2 0.2}{3.8 0.8}
\arr{5.2 0.2}{5.8 0.8}

\arr{1.2 2.2}{1.8 2.8}
\arr{3.2 2.2}{3.8 2.8}
\arr{5.2 2.2}{5.8 2.8}

\arr{0.2 2.8}{0.8 2.2}
\arr{2.2 2.8}{2.8 2.2}
\arr{4.2 2.8}{4.8 2.2}
\arr{6.2 2.8}{6.8 2.2}

\arr{0.2 0.8}{0.8 0.2}
\arr{2.2 0.8}{2.8 0.2}
\arr{4.2 0.8}{4.8 0.2}

\arr{1.2 1.8}{1.8 1.2}
\arr{3.2 1.8}{3.8 1.2}
\arr{5.2 1.8}{5.8 1.2}
\arr{7.2 1.8}{7.8 1.2}

\arr{2.2 0.97}{2.8 0.83}
\arr{4.2 0.97}{4.8 0.83}
\arr{6.2 0.97}{6.8 0.83}

\arr{3.2 0.83}{3.8 0.97}
\arr{5.2 0.83}{5.8 0.97}
\arr{7.2 0.83}{7.8 0.97}
\put{$\Gamma(\Lambda)$} at 9 2

\endpicture} at 5 0
\put{\beginpicture
\setcoordinatesystem units <.5cm,.5cm>
\multiput{$\bullet$} at 0 1  0 3  1 0  1 2  2 1  3 1 /
\plot 2 1  1 0  0 1  1 2 /
\plot 0 3  2 1  3 1 /
\circulararc 360 degrees from 2.22 1 center at 2 1 
\put{$\Cal L$} at -1 2
\endpicture} at 2.8 -2.5
\put{\beginpicture
\setcoordinatesystem units <.5cm,.5cm>
\multiput{$\bullet$} at -1 1  0 1  0 3  1 0  1 2  2 1  /
\plot 2 1  1 0  0 1  1 2 /
\plot 0 3  2 1  /
\plot -1 1  0 1 /
\circulararc 360 degrees from 0.22 1 center at 0 1 
\put{$\Cal K$} at 2.7 2
\endpicture} at 6.4 -2.5
\endpicture}
$$
	\medskip
(3) The local algebra $\Lambda$ with generators $x,y$ and relations $x^2 = y^2$ and
$xy = 0$.  In order to present $\Lambda$-modules, we use here 
the following convention: the bullets represent base vectors, the lines
marked by $x$ or $y$
show that the multiplication by $x$ or $y$, respectively, 
sends the upper base vector to the lower one (all other multiplications by $x$ or
$y$ are supposed to be zero). The upper line shows all the indecomposable modules
in $\Cal L$, the lower one those in $\Cal K.$
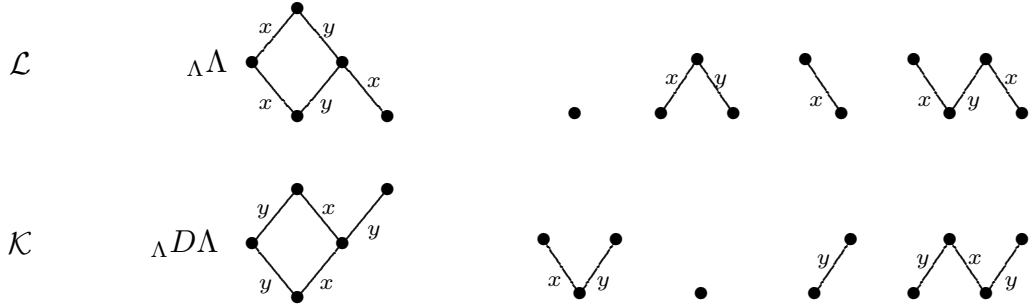
$$
{\beginpicture
\setcoordinatesystem units <1cm,1cm>

\put{\beginpicture
\setcoordinatesystem units <.5cm,.6cm>
\multiput{$\bullet$} at 1 2  0 1  2 1  1 0  3 0 /
\plot 2 1  1 0  0 1  1 2  3 0 /
\put{$\ssize x\strut$} at 0.3 1.7 
\put{$\ssize x\strut$} at 0.3 0.3 
\put{$\ssize y\strut$} at 1.7 1.7 
\put{$\ssize y\strut$} at 1.65 0.3 
\put{$\ssize x\strut$} at 2.7 0.7 
\put{${}_\Lambda\Lambda$} at -1 1

\endpicture} at 0 2

\put{\beginpicture
\setcoordinatesystem units <.4cm,.6cm>
\multiput{$\bullet$} at  0  0 /
\put{} at 1 2 
\endpicture} at 3.3 2

\put{\beginpicture
\setcoordinatesystem units <.4cm,.6cm>
\multiput{$\bullet$} at  0 0  1 1  2 0 /
\plot 0 0  1 1  2 0  /
\put{$\ssize x\strut$} at 0.3 0.7 
\put{$\ssize y\strut$} at 1.65 0.7 
\put{} at 1 2 
\endpicture} at 4.5 2

\put{\beginpicture
\setcoordinatesystem units <.4cm,.6cm>
\multiput{$\bullet$} at  0 1  1 0  /
\plot 0 1 1 0  /
\put{$\ssize x\strut$} at 0.3 0.3 
\put{} at 1 2 
\endpicture} at 5.9 2

\put{\beginpicture
\setcoordinatesystem units <.4cm,.6cm>
\multiput{$\bullet$} at  0 1  2 1  1 0  3 0 /
\plot 3 0  2 1  1 0  0 1  /
\put{$\ssize x\strut$} at 0.3 0.3 
\put{$\ssize x\strut$} at 2.7 0.7 
\put{$\ssize y\strut$} at 1.65 0.3 
\put{} at 1 2 
\endpicture} at 7.5 2

\put{\beginpicture
\setcoordinatesystem units <.5cm,.6cm>
\multiput{$\bullet$} at  1 2  0 1  2 1  1 0  3 2  /
\plot 2 1  1 0  0 1  1 2  2 1  /
\plot 2 1  3 2  /
\put{$\ssize y\strut$} at 0.25 1.7 
\put{$\ssize y\strut$} at 0.3 0.3 
\put{$\ssize x\strut$} at 1.7 1.7 
\put{$\ssize x\strut$} at 1.65 0.3 
\put{$\ssize y\strut$} at 2.7 1.3
\put{${}_\Lambda D\Lambda$} at -1.5 1
\endpicture} at -0.2 0

\put{\beginpicture
\setcoordinatesystem units <.4cm,.6cm>
\multiput{$\bullet$} at  0 1  1 0  2 1 /
\plot 0 1  1 0  2 1  /
\put{$\ssize x\strut$} at 0.3 0.3 
\put{$\ssize y\strut$} at 1.65 0.3 
\put{} at 1 2 
\endpicture} at 3.2 0

\put{\beginpicture
\setcoordinatesystem units <.4cm,.6cm>
\multiput{$\bullet$} at  0  0 /
\put{} at 1 2 
\endpicture} at 4.7 0

\put{\beginpicture
\setcoordinatesystem units <.4cm,.6cm>
\multiput{$\bullet$} at  0 0  1 1  /
\plot 0 0  1 1 /
\put{$\ssize y\strut$} at 0.25 0.7 
\put{} at 1 2 
\endpicture} at 6 0

\put{\beginpicture
\setcoordinatesystem units <.4cm,.6cm>
\multiput{$\bullet$} at  0 0  1 1  2 0  3 1 /
\plot 0 0  1 1  2 0  3 1  /
\put{$\ssize y\strut$} at 0.25 0.7 
\put{$\ssize x\strut$} at 1.7 0.7 
\put{$\ssize y\strut$} at 2.7 0.3 
\put{} at 1 2 
\endpicture} at 7.5 0

\put{$\Cal L$} at -3 2
\put{$\Cal K$} at -3 0

\endpicture}
$$
Let us stress the following: All the indecomposable modules in
$\Cal L\setminus\Cal P$ as well as those in $\Cal K \setminus \Cal Q$ are
$\Lambda'$-modules, where $\Lambda' = k[x,y]/\langle x,y\rangle^2.$ Note that
the category of $\Lambda'$-modules is stably equivalent to the category of
Kronecker modules, thus all its regular components are homogeneous tubes.
In $\Cal L$ we find two indecomposable modules which belong to one tube, in
$\Cal K$ we find two indecomposable modules which belong to another tube.
The algebra $\Lambda'$ has an automorphism which exchanges these two tubes;
this is an outer automorphism, and it {\bf cannot} be lifted to
an automorphism of $\Lambda$ itself. 
	\medskip
(4) In the last example to be presented here, $\Cal L$ (and therefore
also $\Cal K$) will be infinite. We consider the quiver
$$
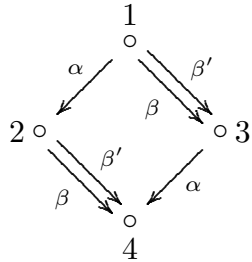
{\beginpicture
\setcoordinatesystem units <1cm,1cm>
\multiput{$\circ$} at 0 1  1 0  1 2  2 1 /
\arr{0.8 1.8}{0.2 1.2}
\arr{1.1 1.8}{1.8 1.1}
\arr{1.2 1.9}{1.9 1.2}
\arr{0.1 0.8}{0.8 0.1}
\arr{0.2 0.9}{0.9 0.2}
\arr{1.8 0.8}{1.2 0.2}
\put{$1$} at 1 2.3
\put{$2$} at -.25 1
\put{$3$} at 2.25 1
\put{$4$} at 1 -.3

\put{$\ssize\alpha$} at 0.4 1.7
\put{$\ssize\beta$} at  1.25 1.25
\put{$\ssize\beta'$} at 1.8 1.7

\put{$\ssize\beta$} at 0.25 0.25
\put{$\ssize\beta'$} at 0.8 0.7
\put{$\ssize\alpha$} at 1.7 0.4

\endpicture}
$$
with the relations $\alpha\beta = \beta\alpha$ and 
$\alpha\beta' = \beta'\alpha$, thus we deal with the tensor product
$\Lambda$ 
of the Kronecker algebra and the path algebra of the quiver of type $\Bbb A_2$
(note that 
tensor products of algebras will be discussed in the second part of this paper
in more detail). For any vertex $i$, we denote by $S(i), P(i), Q(i)$ the simple,
or indecomposable projective or indecomposable injective $\Lambda$-module
corresponding to $i$, respectively. 
	\medskip
The categories $\Cal L$ and $\Cal K$ can be described very 
well using the category of Kronecker modules.
By definition, the Kronecker quiver $K$ has two vertices, a source 
and a sink, and two arrows going from the source to the sink.
Thus a Kronecker module is a quadruple $(U,V,w, w')$
consisting of two vector spaces $U,V$ and two linear maps
$w,w'\:U\to V$.  
We define functors
$\eta,\eta'\:\mod kK \to \mod\Lambda$, sending $M = (U,V,w, w')$ to
the representations
$$
{\beginpicture
\setcoordinatesystem units <1.2cm,1.2cm>
\put{\beginpicture
\arr{0.8 1.8}{0.2 1.2}
\arr{1.1 1.8}{1.75 1.15}
\arr{1.2 1.9}{1.9 1.2}
\arr{0.1 0.8}{0.75 0.15}
\arr{0.2 0.9}{0.9 0.2}
\arr{1.8 0.8}{1.2 0.2}
\put{$0$} at .98 2
\put{$V$} at -.1 1
\put{$U$} at 2.05 .95
\put{$V\oplus V$} at 1 -.1

\put{$\left[\smallmatrix1\cr0\endsmallmatrix\right]$} at 0.1 0.25
\put{$\left[\smallmatrix0\cr1\endsmallmatrix\right]$} at 0.8 0.8
\put{$\left[\smallmatrix w&w'\endsmallmatrix\right]$} at 2 0.4
\put{$\eta(M) =$} at -1.2 1
\endpicture} at 0 0 
\put{\beginpicture
\arr{0.8 1.8}{0.2 1.2}
\arr{1.1 1.8}{1.75 1.15}
\arr{1.3 1.8}{1.9 1.2}
\arr{0.1 0.8}{0.75 0.15}
\arr{0.2 0.9}{0.9 0.2}
\arr{1.8 0.8}{1.2 0.2}
\put{$U\oplus U$} at .98 2
\put{$V$} at -.1 1
\put{$U$} at 2.05 .95
\put{$0$} at 1 -.1

\put{$\left[\smallmatrix1 & 0\endsmallmatrix\right]$} at 1.2 1.25
\put{$\left[\smallmatrix0 & 1\endsmallmatrix\right]$} at 1.85 1.7
\put{$\left[\smallmatrix w\cr w'\endsmallmatrix\right]$} at 0.1 1.7
\put{$\eta'(M) =$} at -1.2 1
\endpicture} at 5 0 
\endpicture}
$$
these functors $\eta, \eta'$ are full embeddings. 

Let us denote by $I$ the indecomposable
injective Kronecker module of length 3, by $T$ the indecomposable projective
Kronecker module of length 3, 
then clearly 
$$
 \eta(I) = \rad P(1) \qquad\text{and} \qquad \eta'(T) = Q(4)/\soc,
$$
and the dimension vector of $\eta(I)$ is $\vect0122$, that of $\eta'(T)$ is 
$\vect2210.$
If $M$ is an indecomposable Kronecker module, then
either $M$ is simple injective and $\eta(M) = S(3),$ or else $M$ is
cogenerated by $I$, and $\eta(M)$ is cogenerated by $\rad P(1),$
thus $\eta(M)$ is a torsionless $\Lambda$-module.
Similarly, either $M$ is simple projective and $\eta'(M) = S(2),$
or else $M$ is generated by $T$ and $\eta'(M)$ is generated by $Q(4)/\soc$,
so that $\eta'(M)$ is divisible.  

On the other hand, nearly all indecomposable
torsionless $\Lambda$-modules are in the image of the functor $\eta$,
the only exceptions are the indecomposable projective modules $P(1), P(3), P(4)$.
Similarly, nearly all indecomposable
divisible $\Lambda$-modules are in the image of the functor $\eta'$,
the only exceptions are the indecomposable injective modules $Q(1), Q(2), Q(4)$.

Altogether, one sees that the category $\Cal L$ has the following
Auslander-Reiten quiver
$$
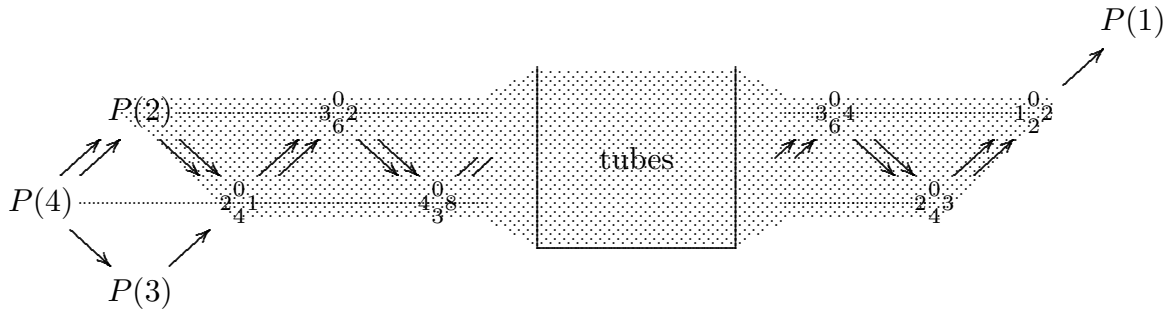
{\beginpicture
\setcoordinatesystem units <1.1cm,1cm>
\put{$P(4)$} at 0 0 
\put{$P(2)$} at 1 1 
\put{$P(3)$} at 1 -1 
\put{$P(1)$} at 11 2 
\put{$\vect0214$} at  2 0
\put{$\vect0326$} at  3 1
\put{$\vect0483$} at  4 0 
\put{$\vect0346$} at  8 1 
\put{$\vect0234$} at  9 0  
\put{$\vect0122$} at 10 1

\arr{0.2 0.3}{0.6 0.7}
\arr{0.4 0.3}{0.8 0.7}
\arr{0.3 -.3}{0.7 -.7}
\arr{1.2 .7}{1.6 0.3}
\arr{1.4 .7}{1.8 0.3}
\arr{1.3 -.7}{1.7 -.3}
\arr{2.2 0.3}{2.6 0.7}
\arr{2.4 0.3}{2.8 0.7}
\arr{3.2 0.7}{3.6 0.3}
\arr{3.4 0.7}{3.8 0.3}
\plot 4.2 0.3  4.4 0.5 /
\plot 4.35 0.3  4.55 0.5 /

\arr{7.4 0.5}{7.6 0.7}
\arr{7.6 0.5}{7.8 0.7}
\arr{8.2 0.7}{8.6 0.3}
\arr{8.4 0.7}{8.8 0.3}
\arr{9.2 0.3}{9.6 0.7}
\arr{9.4 0.3}{9.8 0.7}
\arr{10.3 1.3}{10.7 1.7}
\put{tubes} at 6 0.5
\plot 5 1.5  5 -0.5  7 -.5  7 1.5 /
\setdots <.5mm>
\plot 1.4 1  4.5 1 /
\plot 0.4 0  4.5  0 /
\plot 7.5 1  10 1 /
\plot 7.5  0  9 0 /

\setshadegrid span <.5mm>
\vshade .8 1 1.2 <z,z,,> 2 -.2 1.2 <z,z,,> 4.5 -.2 1.2  <z,z,,>
 5 -.5 1.5 <z,z,,> 7 -.5 1.5 <z,z,,>
 7.5 -.2 1.2  <z,z,,> 9 -.2 1.2 <z,z,,> 10.2 1 1.2 /
\put{} at 12.3 0
\endpicture}
$$
where the dotted part are the torsionless modules which are in the image
of the functor $\eta$.
The category $\Cal L/\Cal P$ is equivalent under $\eta$ to the
category of Kronecker modules without simple direct summands.
	\medskip
Dually, the category $\Cal K$ has the following Auslander-Reiten
quiver:
$$
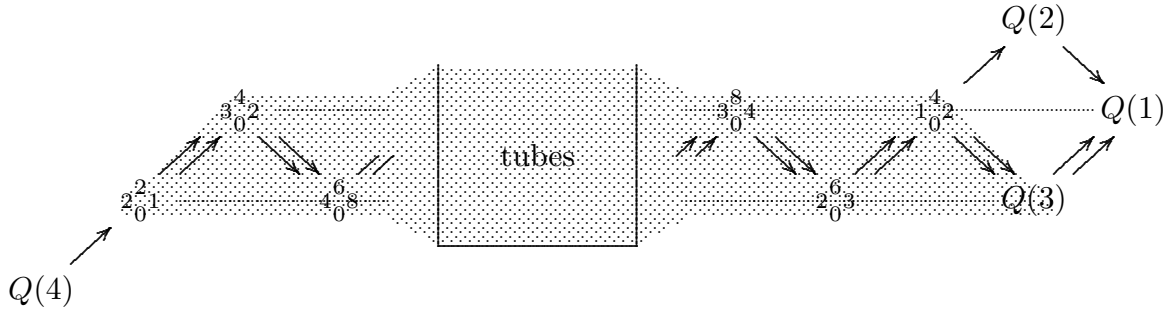
{\beginpicture
\setcoordinatesystem units <1.1cm,1cm>
\put{$Q(4)$} at 1 -1 
\put{$Q(2)$} at 11 2 
\put{$Q(3)$} at 11 0 
\put{$Q(1)$} at 12 1 

\put{$\vect2210$} at  2 0
\put{$\vect4320$} at  3 1
\put{$\vect6480$} at  4 0 
\put{$\vect8340$} at  8 1 
\put{$\vect6230$} at  9 0  
\put{$\vect4120$} at 10 1
\arr{1.3 -.7}{1.7 -.3}
\arr{2.2 0.3}{2.6 0.7}
\arr{2.4 0.3}{2.8 0.7}
\arr{3.2 0.7}{3.6 0.3}
\arr{3.4 0.7}{3.8 0.3}
\plot 4.2 0.3  4.4 0.5 /
\plot 4.35 0.3  4.55 0.5 /

\arr{7.4 0.5}{7.6 0.7}
\arr{7.6 0.5}{7.8 0.7}

\arr{8.2 0.7}{8.6 0.3}
\arr{8.4 0.7}{8.8 0.3}
\arr{9.2 0.3}{9.6 0.7}
\arr{9.4 0.3}{9.8 0.7}
\arr{10.2 0.7}{10.6 0.3}
\arr{10.4 0.7}{10.8 0.3}
\arr{11.2 0.3}{11.6 0.7}
\arr{11.4 0.3}{11.8 0.7}
\arr{10.3 1.3}{10.7 1.7}
\arr{11.3 1.7}{11.7 1.3}
\put{tubes} at 6 0.5
\plot 5 1.5  5 -0.5  7 -.5  7 1.5 /
\setdots <.5mm>
\plot 3.4 1  4.5 1 /
\plot 2.4 0  4.5  0 /
\plot 7.5 1  11.6 1 /
\plot 7.5  0  10.6 0 /

\setshadegrid span <.5mm>
\vshade 1.8 -.2 0 <z,z,,> 3 -.2 1.2 <z,z,,> 4.5 -.2 1.2  <z,z,,>
 5 -.5 1.5 <z,z,,> 7 -.5 1.5 <z,z,,>
 7.5 -.2 1.2  <z,z,,> 10 -.2 1.2 <z,z,,> 11.2 -.2 0 /
\put{} at -.3 0
\endpicture}
$$
here, the dotted part are the divisible modules which are in the image
of the functor $\eta'$ and we see that now the functor $\eta'$ furnishes
an equivalence between the category $\Cal K/\Cal Q$ and the 
category of Kronecker modules without simple direct summands.
	\bigskip
Let us add an interesting property of the functor $\gamma$.
	\medskip
{\bf (3.4)} {\bf Proposition.} {\it Let $M$ be indecomposable, torsionless, but not projective. Then
$\top M$ and $\soc \gamma M$ are isomorphic.}
	\medskip
Proof: In order to calculate $\gamma M = \Sigma\tau M$, we start with a minimal
projective presentation $f\:P_1 \to P_0$, apply the functor $\nu =D\Hom(-,\Lambda)$
to $f$ and take as $\gamma M$ the image of $\nu(f).$ Here, the embedding of $\gamma M$
into $\nu P_0$ is an injective envelope. Since $P_0$ is a projective cover of $M$,
we have $\top P_0 \simeq \top M$; since $\nu P_0$ is an injective envelope of $\gamma M$,
we have $\soc \gamma M \simeq \soc \nu P_0$. And of course, we have $\top P_0 \simeq
\soc \nu P_0.$
	\medskip
This property of $\gamma$ is nicely seen in the last example!
Of course,
the canonical bijection between the indecomposable projective and the
indecomposable injective modules has also this property. 
	\bigskip\bigskip
{\bf 4. The representation dimension of a torsionless-finite artin algebra} 
	\medskip

{\bf (4.1) Theorem.} {\it Let $\Lambda$ be a torsionless-finite artin algebra. Let $M$ be the direct sum
of all indecomposable $\Lambda$-modules which are torsionless or
divisible, one from each isomorphism class. Then the global dimension of
$\End(M)$ is at most $3$.}
	\medskip
Note that such a module $M$ is a generator-cogenerator, thus we see:
{\it If $\Lambda$ is in torsionless-finite and
representation-infinite, then the direct sum of all $\Lambda$-modules
which are torsionless or divisible is an Auslander generator.} In particular: 
	\medskip
{\bf (4.2) Corollary.} {\it  If $\Lambda$ is a torsionless-finite artin algebra, then
$\repdim\Lambda \le 3$.}
	\bigskip

For the proof of Theorem 4.1, we need the following lemma which  
goes back to Auslander's Queen Mary notes [A] where is was
used implicitly. The formulation is due to [EHIS] and [CP].
Given modules $M,X$, denote by $\Omega_M(X)$ the kernel of a minimal right $\add M$-approximation 
$g_{MX}\:M'\to X$. By definition, the $M$-dimension $M$-$\dim X$ 
of $X$ is the minimal value $i$ such that $\Omega_M^i(X)$ belongs to $\add M$.
	\medskip
{\bf (4.3) Auslander-Lemma.} {\it Let $M$ be a $\Lambda$-module. If 
$M$-$\dim X \le d$ for all $\Lambda$-modules $X$, then the global
dimension of $\End(M)$ is less or equal $d+2$.} If $M$ is a 
generator-cogenerator, then also the converse holds: if the global dimension of $\End(M)$
is less or equal $d+2$ with $d\ge 0$, then $M$-$\dim X \le d$ for all $\Lambda$-modules $X$.
	\medskip
Let us outline the proof of the first implication. 
Thus, let us assume $M$-$\dim X \le d$ for all $\Lambda$-modules $X$.
Given any $\End(M)$-module $Y$, we want to construct
a projective $\End(M)$-resolution of length at most $d+2$. 
The projective $\End(M)$-modules are of the form $\Hom(M,M')$ with $M'\in \add M.$
Consider a projective presentation of $Y$, thus an exact sequence
$$
 \Hom(M,M'') @>\phi>> \Hom(M,M') \to Y \to 0,
$$
with $M',M'' \in \add M$. Note that $\phi = \Hom(M,f)$ for some map $f\:M'' \to M'$.
Let $X$ be the kernel of $f$, thus $\Hom(M,X)$ is the kernel of $\Hom(M,f) = \phi.$
Inductively, we construct minimal right $\add M$-approximations
$$
 M_{i} @>g_i>> \Omega_M^i(X),
$$
starting with $\Omega_M^0(X) = X$, so that the kernel of $g_i$ is just 
$\Omega_M^{i+1}$, say
with inclusion map $u_{i+1}\:\Omega_M^{i+1} \to M_i.$
Thus, we get a sequence of maps
$$
 0 \to \Omega_M^d(X) @>u_d>> M_{d-1} @>u_{d-1}g_{d-1}>> M_{d-2} @>>> \cdots @>>> M_1 @>u_1g_1>> M_0 @>f_0>> X @>>> 0 
$$
If we apply the functor $\Hom(M,-)$ to this sequence, we get 
an exact sequence
$$
 0 \to \Hom(M,\Omega_M^d(X)) \to \Hom(M,M_{d-1}) \to \cdots \to \Hom(M,M_0) \to \Hom(M,X) \to 0
$$
(here we use that we deal with right $M$-approximations and that $\Hom(M,-)$ is
left exact). Since we assume that $\Omega_M^d(X)$ is in $\add M$, we see that we
have constructed a projective resolution of $\Hom(M,X)$ of length $d$. 
Combining this with the exact sequence
$$
  0 \to \Hom(M,X) \to  \Hom(M,M'') @>\phi>> \Hom(M,M') \to Y \to 0,
$$ 
we obtain a projective resolution of $Y$ of length $d+2$. This completes the proof.
	\medskip
	
{\bf (4.4) Proof of Theorem 4.1.} As before, let 
$\Cal L$ be the class of torsionless $\Lambda$-modules, and
$\Cal K$ be the class of divisible $\Lambda$-modules. Since $\Lambda$
is torsionless-finite there are $\Lambda$-modules $K,L$ 
with $\add K = \Cal K$, and $\add L = \Cal L$. 
Let $M = K\oplus L.$ We use the Auslander Lemma.

Let $X$ be a $\Lambda$-module. 
Let $U$ be the trace of $\Cal K$ in $X$ (this is the sum
of the images of maps $K \to X$). Since $\Cal K$ is closed under direct sums and factor modules,
$U$ belongs to $\Cal K$ (it is the largest submodule of $X$ which belongs to $\Cal K$).
Let $p\:V \to X$ be a right $\Cal L$-approximation of $X$ (it exists, since we assume that
$\Cal L$ is finite). Since $\Cal L$ contains all the projective modules, it follows
that $p$ is surjective. Now we form the pullback
$$
\CD
V @>p>> X \cr
@Au'AA @AAuA \cr
W @>>p'> U 
\endCD
$$
where $u\:U \to X$ is the inclusion map.
With $u$ also $u'$ is injective, thus $W$ is a submodule of $V\in \Cal L$. Since $\Cal L$
is closed under submodules, we see that $W$ belongs to $\Cal L$.
On the other hand, the pullback gives rise to the exact sequence
$$
 0 @>>> W @>{\bmatrix p' &\! -u'\endbmatrix}>> U\oplus V @>{\bmatrix u \cr p\endbmatrix}>> X @>>> 0
$$
(the right exactness is due to the fact that $p$ is surjective). By construction, the map $\bmatrix u \cr p\endbmatrix$
is a right $M$-approximation, thus $\Omega_{M}(X)$ is a direct summand
of $W$ and therefore in $\Cal L \subseteq \add M.$ This completes the proof.

	\bigskip\bigskip
\vfill\eject
{\bf 5. Classes of torsionless-finite artin algebras} 
	\medskip
In the following, let $\Lambda$ be an artin algebra with radical $J$.
	\medskip
Before we deal with specific classes of torsionless-finite artin algebras,
let us mention two characterizations of torsionless-finite artin algebras:
	\medskip
{\bf (5.1) Proposition.} {\it An artin algebra $\Lambda$ is torsionless-finite
if and only if there exists a faithful module 
$M$ such that the subcategory of modules cogenerated by $M$ is finite.}
	\medskip
Proof. If $\Lambda$ is torsionless-finite, we can take 
$M = {}_\Lambda\Lambda$. Conversely, assume that $M$ is faithful and that
the subcategory of modules cogenerated by $M$ is finite.
Since $M$ is faithful, the regular representation ${}_\Lambda\Lambda$ itself
is cogenerated by $M$, thus all the torsionless-finite $\Lambda$-modules are 
cogenerated by $M$. This shows that $\Lambda$ is torsionless-finite.
	\medskip
Actually, also non-faithful modules can similarly be used in order to 
characterize torsionless-finiteness, for example we can take $\rad\Lambda$ (note
that for a non-zero artin algebra, $\rad\Lambda$ is never faithful,
since it is annihilated by the right socle of $\Lambda$):
	\medskip
{\bf (5.2) Proposition.} {\it An artin algebra $\Lambda$ is torsionless-finite
if and only if 
the subcategory of modules cogenerated by $J$ is finite.}
	\medskip
Proof: On the one hand, modules cogenerated by $\rad\Lambda$ are torsionless. 
Conversely, assume that there are only finitely many isomorphism classes of 
indecomposable $\Lambda$-modules which are cogenerated by $\rad\Lambda$. Then
$\Lambda$ is torsionless-finite, according to following lemma.	\bigskip

{\bf (5.3) Lemma.} {\it Let $N$ be an indecomposable torsionless $\Lambda$-module. 
Then either $N$ is projective or else $N$ is cogenerated by $J.$}
	\medskip
Proof: 
Let $N$ be indecomposable and torsionless, but not cogenerated by
$J$. We claim that $N$ is projective. 
Since $N$ is torsionless, there is an inclusion map
$u\:N \to P = \bigoplus P_i$ with indecomposable projective modules $P_i$.
Let $\pi_i\:P \to P_i$ be the canonical projection onto the direct summand
$P_i$ of $P$ and $\epsilon_i\:P_i\to S_i$ the canonical projection
of $P_i$ onto its top.
If $\epsilon\pi_iu = 0$
for all $i$, then $N$ is contained in the radical of $P$, thus 
cogenerated by $\rad\Lambda$, a contradiction. Thus there is some
index $i$ with $\epsilon\pi_iu \neq 0$, but this implies that
$\pi_i u$ is surjective. Since this is a surjective map onto a
projective module, we see that $\pi_iu$ is a split epimorphism. But
we assume that $N$ is indecomposable, thus $\pi_iu$ is an isomorphism.
This shows that $N \simeq P_i$ is projective. Altogether, we see
that there are only finitely many isomorphism classes of indecomposable
torsionless $\Lambda$-modules, namely those cogenerated by $\rad\Lambda$,
as well as some additional ones which are projective. 
	\bigskip
Here are now some classes of torsionless-finite artin algebras:
	\medskip
{\bf (5.4) Artin algebras $\Lambda$ with 
$\Lambda/\soc(\Lambda_\Lambda)$ representation-finite.} 
Let $N$ be an indecomposable torsionless $\Lambda$-module which is not projective.
By Lemma 5.3, there is an embedding $u\:N \to J^t$ for some $t$. 
Let $I = \soc(\Lambda_\Lambda).$ Then $u(IN) = Iu(N) \subseteq I(J^t) = 0$,
thus $IN = 0.$ This shows that $N$ is a $\Lambda/I$-module. Thus $N$ belongs to
one of the finitely many isomorphism classes of indecomposable $\Lambda/U$-modules.
This shows that $\Lambda$ is torsionless-finite. 

{\it If $J^n = 0$ and $\Lambda/J^{n-1}$ is representation-finite, then 
$\Lambda$ is torsionless-finite.} Namely, $J^{n-1} \subseteq \soc(\Lambda_\Lambda),$
thus, if $\Lambda/J^{n-1}$ is representation-finite, also its factor algebra 
$\Lambda/\soc(\Lambda_\Lambda)$ is torsionless-finite.
This shows: 
{\it If $J^n = 0$ and $\Lambda/J^{n-1}$ is representation-finite, then 
the representation dimension of
$\Lambda$ is at most $3$.} (Auslander [A], Proposition, p.143)
	\medskip
{\bf (5.5) Artin algebras with radical square zero.} Following Auslander (again [A], Proposition, p.143)
This is the special case $J^2 = 0$ of 5.4. Of course, here the proof of the
torsionless-finiteness is very easy: An indecomposable torsionless module is either projective
or simple. Similarly, an indecomposable divisible module is either injective or
simple, and any simple module is either torsionless or divisible. Thus the
module $M$ exhibited in Theorem 4.1 is the direct sum of all indecomposable projective, 
all indecomposable injective, and all simple modules. 
	\medskip
{\bf (5.6) Minimal representation-infinite algebras.} Another special case of 5.4 is of interest:
We say that $\Lambda$ is minimal representation-infinite provided
$\Lambda$ is representation-infinite, but any proper factor algebra is representation-finite. If $\Lambda$ is minimal representation-infinite, and $n$ is minimal
with $J^n = 0,$ then $\Lambda/J^{n-1}$ is a proper factor algebra, thus
representation-finite. 
	\medskip
{\bf (5.7) Hereditary artin algebras.} If $\Lambda$ is hereditary, then the only torsionless modules are the projective
modules and the only divisible modules are the injective ones, thus 
the module $M$ of Theorem 4.1 is the direct sum of all indecomposable modules
which are projective or injective. 
In this way, we recover Auslander's result 
([A], Proposition, p. 147). 
	\medskip
{\bf (5.8) Artin algebras stably equivalent to hereditary algebras.} Let $\Lambda$ be stably 
equivalent to a hereditary artin algebra. Then 
an indecomposable torsionless module is either
projective or simple ([AR1], Theorem 2.1), thus there are only finitely many
isomorphism classes of torsionless $\Lambda$-modules.
Dually, an indecomposable divisible module is either
injective or simple.
Thus, again we see the structure of the module $M$ of Theorem 4.1
and we recover Proposition 4.7
of Auslander-Reiten [AR2].
	\medskip
{\bf (5.9) Right glued algebras} (and similarly left glued algebras): An artin
algebra $\Lambda$ is said to be {\it right glued}, provided the functor
$\Hom(D\Lambda,-)$ is of finite length, or equivalently, provided almost all
indecomposable modules have projective dimension equal to 1.
The condition that $\Hom(D\Lambda,-)$ is of finite length implies that there
are only finitely many isomorphism classes of divisible $\Lambda$-modules.
Also, the finiteness of the isomorphism classes
of indecomposable modules of projective dimension greater than 1 implies 
torsionless-finiteness. We see that {\it right glued algebras have representation
dimension at most $3$} (a result of Coelho-Platzeck [CP]).
	\medskip
{\bf (5.10) Special biserial algebras without indecomposable projective-injective modules.} 
In order to show that these artin algebras are torsionless-finite, we 
need the following Lemma.
	\medskip
{\bf Lemma.} {\it Let $\Lambda$ be special biserial and $M$ a $\Lambda$-module.
The following assertions are equivalent.
\item{\rm (i)} $M$ is a direct sum of local string modules.
\item{\rm (ii)} $\alpha M \cap \beta M = 0$ for arrows $\alpha \neq \beta.$}
	\medskip
Proof: (i) $\implies$ (ii). We can assume that $M$ is indecomposable, but then it 
should condition (ii) is satisfied.

(ii) $\implies$ (i): We can assume that $M$ is indecomposable. For a band module,
condition (ii) is clearly not satisfied. And for a string module $M$, condition (ii) is 
only satisfied in case $M$ is local.
	\bigskip
Proof that special biserial algebras without indecomposable projective-injective modules
are torsionless-finite: Assume that $\Lambda$ is special biserial and that there is
no indecomposable projective-injective module.
Then all the indecomposable projective modules are string modules (and of course
local). Thus any projective module satisfies the condition (i) and therefore also 
the condition (ii). But if a module $M$ satisfies the condition (ii), also every
submodule of $M$ has this property. This shows that all torsionless
modules satisfy the condition (ii). It follows that indecomposable torsionless modules are
local string modules, and the number of such modules is finite. 
	\bigskip

It follows from 5.10 that {\it all special biserial algebras 
have representation dimension at most $3$,} as shown in [EHIS]. For the proof 
one uses the following general observation (due to [EHIS] 
in case the representation dimension is $3$): 
	\medskip
{\bf (5.11) Proposition.} {\it Let $\Lambda$ be an artin algebra. Let
$P$ be indecomposable projective-injective $\Lambda$-module. There is 
a minimal two-sided ideal $I$ such that $IP \neq 0.$
Let $\Lambda' = \Lambda/I.$ 
Then either $\Lambda'$ is semisimple or else $\repdim \Lambda \le \repdim \Lambda'.$}
	\medskip
Proof: Note that all the indecomposable $\Lambda$-modules not isomorphic to $P$
are annihilated by $I$, thus they are $\Lambda'$-modules.

First assume that $\Lambda$ is representation finite, 
thus $\repdim \Lambda \le 2.$ Now $\Lambda'$ is also representation finite, and by assumption not semisimple, thus $\repdim \Lambda' = 2$.  This yields the claim. (Actually, $\Lambda$ cannot be semisimple, since otherwise also $\Lambda'$
semisimple, thus $\repdim \Lambda = 2$ and therefore $\repdim \Lambda = \repdim \Lambda'.$)

Now assume that $\Lambda$ is not representation finite, with representation dimension $d$.
Let $M'$ be an Auslander generator for $\Lambda'$, thus, according the second assertion
of the Auslander-Lemma asserts that $M'$-$\dim X \le d$ for all $\Lambda'$-modules $X$.
Let $M = M'\oplus P$. This is clearly a generator-cogenerator. 
We want to show the any indecomposable $\Lambda$-module has
$M$-dimension at most $d$ (then $\End(M)$ has global dimension at most $d$
and therefore the representation dimension of $\Lambda$ is at most $d$).

Let $X$ be an indecomposable $\Lambda$-module.
Now $X$ may be isomorphic to $P$, then $X$ is in $\add M$,
thus its $M$-dimension is $0$. 

So let us assume that $X$ is not isomorphic to
$P$, thus a $\Lambda'$-module. Let $g\:M'' \to X$ be a minimal right $M'$-approximation of
$X$. We claim that {\it $g$ is even a minimal right $M$-approximation.} Now $M''$ is
in $\add M$, thus we only have to show that any map $f\:M_i \to X$ factors through $g$,
where $M_i$ is an indecomposable direct summand of $M$. This is clear in case $M_i$ 
is a direct summand of $M'$, thus we only have to look at the case $M_i = P$. But
since $X$ is a annihilated by $I$, the map $f\:P \to X$ vanishes on $IP$, thus $f$
factors through the projection map $p\:P \to P/IP$, say $f = f'p$ with
$f'\:P/IP \to X.$ Since $P/IP$ is an indecomposable projective $\Lambda'$-module, 
it belongs to $\add M'$, thus $f'$ factors through $g$, say $f' = gf''$
for some $f''\:P/IP \to M''.$ Thus $f = f'p = gf''p$ factors through $g$.
This concludes the proof that $g$ is a minimal right $M$-approximation.

Now $\Omega_M(X)$ is the kernel of $g$, thus $\Omega_M(X) = \Omega_{M'}(X)$,
in particular, this is again
a $\Lambda'$-module. Thus, inductively we see that $\Omega^i_M(X) = \Omega^i_{M'}(X)$
for all $i$. But we know that $\Omega^d_M(X) = \Omega^d_{M'}(X)$ is in $\add M'$,
and $\add M' \subseteq \add M$. This shows that $X$ has $M$-dimension at most $d$.
	\bigskip
There are many other classes of artin algebras studied in the literature which
can be shown to be torsionless-finite, thus have representation dimension at most $3$
(note that also Theorem 5.1 of [X1] deals with artin algebras which are divisible-finite,
thus torsionless-finite). 
	\bigskip
{\bf (5.12) Further algebras with representation dimension $3$.}
We have seen that many artin algebras of interest are torsionless-finite
and thus their representation dimension is at most $3$.
But we should note that not all artin algebras with representation dimension
at most $3$ are torsionless-finite.

Namely, it is easy to construct special biserial algebras which are not torsionless-finite.
And there are also many tilted algebras as well as canonical algebras
which are not torsionless-finite, whereas all tilted and all 
canonical algebras have representation-dimension at most $3$,
see Assem-Platzeck-Trepode [APT] and Oppermann [O3].
Actually, as Happel-Unger [HU] have shown, all piecewise hereditary
algebras have representation dimension at most $3$ (an algebra $\Lambda$
is said to be piecewise hereditary provided the derived category $D^b(\mod\Lambda)$
is equivalent as a triangulated category to the bounded derived category of some
hereditary abelian category), but of course not all are torsionless-finite.  
	\bigskip\bigskip
\vfill\eject
\centerline{\gross Part II. The Oppermann dimension}
	\smallskip
\centerline{\gross and tensor products of algebras.} 
	\bigskip

We consider now $k$-algebras $\Lambda$, where $k$ is a field. 
	\medskip
{\bf 6. Oppermann dimension.} 
	\medskip

{\bf (6.1)} 
Let $R = k[T_1,\dots,T_d]$ be the polynomial ring in $d$ variables with coefficients
in $k$ and $\Max R$ its maximal spectrum, this is the
set of maximal ideals of $R$  endowed with the Zariski topology.
For example, given $\alpha = (\alpha_1,\dots,\alpha_d)\in k^d$, there is
the maximal ideal 
$\m_\alpha = \langle T_i-\alpha_i\mid 1 \le i \le d\rangle.$ 
In case $k$ is algebraically closed, we may 
identify in this way $k^d$ with $\Max R$, otherwise $k^d$ yields 
only part of $\Max R.$ In general, we will denote an element of $\Max R$ 
by  $\alpha$ (or also by $\m_\alpha$, if we want to stress that we consider $\alpha$
as a maximal ideal), and $R_\alpha$ will 
denote the corresponding localization of $R$, whereas
$S_\alpha = R/\m_\alpha$ is the corresponding simple $R$-module (note that all simple
$R$-modules are obtained in this way).
For any ring $A$, let
$\fin A$ be the category of finite length $A$-modules.

By definition, a {\it $\Lambda\otimes_kR$-lattice} $L$ is a finitely generated $\Lambda\otimes_kR$-module which is projective (thus free) as an $R$-module, we also
will say that $L$ is a {\it $d$-dimensional lattice} for $\Lambda$.
Given a $\Lambda\otimes_kR$-lattice $L$,
we may look at the functor
$$
 L\otimes_R - \: \fin R \longrightarrow \mod\Lambda.
$$
Since $L_R$ is projective, this is an exact functor. This means that given an 
exact sequence of $R$-modules, applying $L\otimes_R -$ we obtain an exact sequence
of $\Lambda$-modules. Thus, if $M,N$ are $R$-modules, and $d$ is a natural number, 
then looking at an element of $\Ext^d_R(M,N)$, we may interprete this element
as the equivalence class $[\epsilon]$ of a long exact sequence $\epsilon$ 
starting with $N$ and ending in $M$, and we may apply $L\otimes_R-$ to $\epsilon$.
We obtain in this way a long exact sequence $L\otimes_R\epsilon$ starting
with $L\otimes_RN$ and ending with $L\otimes_RM$ and its equivalence class
$[L\otimes_R\epsilon]$ in $\Ext^d_\Lambda(L\otimes_RM,L\otimes_RN).$ Since
this equivalence class $[L\otimes_R\epsilon]$ only depends on $[\epsilon]$, we 
obtain the following function, also denoted by $L\otimes_R-$:
$$
  L\otimes_R-\:\Ext^d_R(M,N) \longrightarrow \Ext^d_\Lambda(L\otimes_RM,L\otimes_RN),
  \quad\text{with}\quad (L\otimes_R-)[\epsilon] = [L\otimes_R\epsilon].
$$

We say that $L$ is a
{\it $d$-dimensional 
Oppermann lattice} for $\Lambda$ provided the set of $\alpha \in \Max R$ such that
$$
 (L\otimes_R -)\left(\Ext^d_R(\fin R_\alpha, \fin R_\alpha)\right) \neq 0,
$$
is dense in $\Max R$; this means that for these $\alpha \in \Max R$,
there are modules $M,N \in \fin R_\alpha$
with  
$$
 (L\otimes_R -)\left(\Ext^d_R(M,N)\right) \neq 0.
$$

Actually, instead of looking at all the modules in $\fin R_\alpha$, 
it is sufficient to deal with the simple module $S_\alpha$.
One knows
that $\Ext_R^d(S_\alpha,S_\alpha)$ is generated as a $k$-space 
by the equivalence class of a long exact sequence of the form
$$
 \epsilon_\alpha\: 
 \quad 0 \to S_\alpha \to M_1  \to \cdots \to M_d \to S_\alpha \to 0
$$
with $R_\alpha$-modules $M_i$ which are indecomposable and of length 2.
If we tensor this exact sequence with $L$, we obtain an exact sequence
$$
 L\otimes_R\epsilon_\alpha\: 
 \quad 0 \to L\otimes_RS_\alpha \to 
 L\otimes_RM_1  \to \cdots \to L\otimes_RM_d \to L\otimes_RS_\alpha 
 \to 0
$$ 
which yields an equivalence class $[L\otimes_R  \epsilon_\alpha]$ in 
$\Ext_\Lambda^d(L\otimes_R S_\alpha,L\otimes_R S_\alpha).$
It is easy to see that the following conditions are equivalent:
\item{(i)} $(L\otimes_R -)\left(\Ext^d_R(\fin R_\alpha, \fin R_\alpha)\right) \neq 0$,
\item{(ii)} $[L\otimes_R\epsilon_\alpha] \neq 0$ as an element of 
  $\Ext_\Lambda^d(L\otimes_R S_\alpha,L\otimes_R S_\alpha).$
	\smallskip
Thus we see: {\it The $d$-dimensional lattice $L$ is an
Oppermann lattice for $\Lambda$ provided the set of $\alpha \in \Max R$ such that
$[L\otimes_R \epsilon_\alpha]$ is a non-zero element of 
$\Ext^d_R(S_\alpha, S_\alpha)$ is dense in $\Max R$.}
	\medskip
By definition, the {\it Oppermann dimension} $\Odim\Lambda$ 
of $\Lambda$
is the supremum of $d$ such that there exists a $d$-dimensional Oppermann lattice
$L$ for $\Lambda.$ 

	\medskip
{\bf (6.2) Examples.} {\it 
\item{\rm (a)} Let $\Lambda$ be the path algebra of a representation-infinite
quiver. Then $\Odim \Lambda = 1.$
\item{\rm (b)} Let $\Lambda$ be a representation-infinite $k$-algebra, where $k$
is an algebraically closed field. Then $\Odim \Lambda \ge 1.$}
	\medskip
Proof. (a) The usual construction of one-parameter families of indecomposable 
$\Lambda$-modules for a representation-infinite quivers shows that $\Odim\Lambda \ge 1.$
On the other hand, the path algebra of a quiver is hereditary, thus $\Ext^2_\Lambda = 0.$
This shows that the Oppermann dimension can be at most $1$.

(b) This follows from the proof of the second Brauer-Thrall conjecture by Bautista [Bt] 
and [Bo], see also [Bo2] and [R5], [R6].
	\medskip

The following result of Oppermann 
([O1], Corollary 3.8) shows that $\Odim \Lambda$ is always finite and that one
obtains in this way an
interesting lower bound for the representation dimension:
	\medskip
{\bf (6.3) Theorem (Oppermann).} {\it Let $\Lambda$ be a finite-dimensional $k$-algebra which is not semisimple. Then}
$$
 \Odim \Lambda + 2 \le \repdim \Lambda.
$$
	\medskip
One may ask whether one always has the equality
$\Odim \Lambda + 2 = \repdim \Lambda,$
this can be considered as a formidable extension of the assertion of the 
second Brauer-Thrall conjecture.

	\bigskip\bigskip
{\bf 7. Tensor products of artin algebras.} 
	\medskip
Quite a long time ago, Changchang Xi [X1] has shown the
following inequality: Given finite-dimensional $k$-algebras $\Lambda, \Lambda'$,
$$ 
  \repdim \Lambda\otimes_k \Lambda' \le \repdim \Lambda + \repdim\Lambda',
$$
provided $k$ is a perfect field. 
This provides an upper bound for the representation dimension of 
$\Lambda\otimes_k \Lambda'.$ But there is also a lower bound, which uses
the Oppermann dimension. Let us draw attention to the following fact:
	\medskip
{\bf (7.1) Theorem.} {\it Let $\Lambda, \Lambda'$ be finite-dimensional $k$-algebras.
Let $L$ be an Oppermann lattice for $\Lambda$ and $L'$ an Oppermann lattice for
$\Lambda'$. Then $L\otimes_k L'$ is an Oppermann lattice for $\Lambda\otimes_k
\Lambda'$.}
	\medskip
Proof: Theorem 7.1 is an immediate consequence of 
Theorem 3.1 in Chapter XI of Cartan-Eilenberg [CE].

Namely, let $L$ be a $d$-dimensional Oppermann lattice for $\Lambda$ and
$L'$ a $d'$-dimensional Oppermann lattice for $\Lambda'$. 
Thus, $L$ is an $R$-lattice with $R = k[T_1,\dots,T_d]$ and
say $L'$ is an $R'$-lattice, where $R' = k[T'_1,\dots,T'_{d'}]$
(with new variables $T'_i$). For $\alpha\in \Max R$ we choose an
exact sequence $\epsilon_\alpha$ such that its equivalence class 
$[\epsilon_\alpha]$ generates $\Ext_R^d(S_\alpha,S_\alpha)$;
similarly, for 
$\alpha'\in \Max R'$ we choose an
exact sequence $\epsilon_{\alpha'}$ such that its equivalence class 
$[\epsilon_{\alpha'}]$ generates $\Ext_R^d(S_{\alpha'},S_{\alpha'})$;

Since $L$ is an Oppermann lattice  for $\Lambda$,
 the set of 
elements $\alpha\in \Max R$ such that $[L\otimes_R\epsilon_\alpha] \neq 0$ is
dense in $\Max R.$ Similarly, since $L'$ is an Oppermann lattice for $\Lambda'$,
the set of 
elements $\alpha'\in \Max R'$ such that $[L\otimes_R\epsilon_{\alpha'}] \neq 0$ is
dense in $\Max R'.$

Now $L\otimes_k L'$ is a $\Lambda'\otimes_k \Lambda' \otimes_k R\otimes_k R'$-lattice
and we may look at 
$$
 (L\otimes_kL') \otimes_{R\otimes R'} (\epsilon_\alpha\vee \epsilon_{\alpha'}).
$$
We claim that its equivalence class is non-zero in 
$
 \Ext_{\Lambda\otimes\Lambda'}^{d+d'}(S_{(\alpha,\alpha')},S_{(\alpha,\alpha')}).
$
This is a special case of Theorem XI.3.1 of Cartan-Eilenberg which
asserts the
following: Let $\Lambda, \Lambda'$ be left noetherian $k$-algebras, where $k$ is
a semisimple commutative ring. Let $M$ be a finitely generated $\Lambda$-module
and  $M'$
a finitely generated $\Lambda'$-module. Then the canonical map 
$$
 \vee\:\Ext_\Lambda^d(M,N) \otimes_k \Ext_{\Lambda'}^{d'}(M',N') \longrightarrow
 \Ext_{\Lambda\otimes_k\Lambda'}^{d+d'}(M\otimes_k M',N\otimes_k N')
$$
is an isomorphism for any $\Lambda$-module $N$, $\Lambda'$-module $N'$ and all
$d,d'\in \Bbb N$.
	\medskip
It remains to note that for dense subsets $X$ of $\Max R$ and $X'$ of $\Max R'$, the
product $X\times X'$ is of course dense in $\Max R\otimes_kR'.$ 
	\bigskip

{\bf (7.2) Corollary.} {\it Let $\Lambda, \Lambda'$ be finite-dimensional $k$-algebras. Then}
$$
 \Odim \Lambda\otimes_k \Lambda' \ge \Odim \Lambda + \Odim \Lambda'.
$$
	\medskip
Note that it is easy to provide examples where we have strict inequality: just take
representation-finite algebras $\Lambda, \Lambda'$ such that the Oppermann
dimension of $\Lambda\otimes_k\Lambda'$ is at least 1, for example consider
$\Lambda = \Lambda' = k[T]/\langle T^2\rangle,$ or take $\Lambda, \Lambda'$ 
path algebra of quivers of type $\Bbb A_n$ with $n \ge 3.$ Note that the 
representation type of the tensor product of any two nonsimple connected $k$-algebras
with $k$ an algebraically closed field, has been determined
by Leszczy\'nski and A.Skowro\'nski [LS].
	\medskip
The combination of the inequalities 6.3 and 7.2 yields:
	\medskip
{\bf (7.3) Corollary.} {\it Let $\Lambda_1, \dots, \Lambda_n$ be finite-dimensional $k$-algebras.
Then}
$$
 \repdim \Lambda_1\otimes_k \cdots \otimes_k \Lambda_n \ge 
  2 + \sum\nolimits_{i=1}^n\Odim \Lambda_i.
$$
	\medskip
In particular, we see: 
{\it Let $\Lambda$ be the tensor product of $n$ $k$-algebras
with Oppermann dimension greater or equal to $1$. Then $\rep\dim \Lambda \ge n+2$.}
Using 6.2, we see:
	\medskip
{\bf (7.4) Corollary.} 
\item{\rm(a)} {\it If $\Lambda$ is the tensor product of $n$ path algebras
of representation-infinite quivers, then $\repdim\Lambda \ge n+2.$
\item{\rm(b)} If $k$ is an algebraically closed field and $\Lambda$ is the
tensor product of $n$ representation-infinite $k$-algebras, then 
$\repdim\Lambda \ge n+2.$}
	\bigskip
As Happel has pointed out, the following remarkable consequence should be stressed:
If $\Lambda_1, \Lambda_2$ are representation-infinite path algebras, then
$\Lambda_1\otimes\Lambda_2$ is never a tilted algebra. After all, tilted algebras
have representation dimension at most $3$, whereas we have shown that 
$\repdim \Lambda_1\otimes_k\Lambda_2 \ge 4.$
	\bigskip\bigskip
{\bf 8. Nicely  tiered algebras.}
	\medskip
Let $Q$ be a finite connected quiver. We say that $Q$ is {\it tiered with $n+1$ tiers}
provided there is a surjective function $l\:Q_0 \to [0,n] = \{z\in \Bbb Z\mid 0 \le z \le n\}$
such that for any arrow $x\to y$ one has $l(x) = l(y)+1.$ Such a function $l$, if it 
exists, is uniquely determined and is called the {\it tier function} for $Q$ and $l(x)$
is said to be the {\it tier} (or the tier number) of the vertex $x$. 
We say that $Q$ is {\it nicely  tiered with $n+1$ tiers} provided $Q$ is tiered with
$n+1$ tiers, say with tier function $l$ such that 
$l(x) = 0$ for all sinks $x$, and $l(x) = n$
for all sources $x$. Clearly, $Q$ is nicely tiered if and only if 
$Q$ has no oriented cyclic paths and any maximal path has length $n$. The tier function $l$
of a nicely tiered quiver $Q$ can be characterized as follows: 
the tier number
$l(x)$ of a vertex $x$ is the length of any maximal path starting in $x$.

	\medskip
Let $Q$ be a nicely  tiered quiver and $M$ a representation of $Q$.
We denote by $M|[a,b]$ the restriction of $M$ to the subquiver of all vertices $x$ with 
$a \le l(x) \le b$. We say that a module lives in the interval $[a,b]$ provided 
$M = M|[a,b].$ 
We say that a module $M$ is generated in tier $a$ provided its top lives in $[a,a].$
Dually, $M$ is said to be cogenerated in tier $a$ provided its socle lives in $[a,a].$
Given a module $M$, and $t\in \Bbb N_0$, let 
${}_tM$ be its $t$-th socle (thus, there is the sequence of submodules
$$
 0 = {}_0M \subseteq {}_1M \subseteq \cdots \subseteq {}_tM \subseteq M
$$
such that ${}_tM/{}_{t-1}M = \soc M/{}_{t-1}M$ for all $t\ge 1$).
	\medskip
We say that an algebra is {\it nicely  tiered} provided it is
given by a nicely  tiered quiver and a set of commutativity relations. Let
$\Lambda$ be nicely tired. Then 
it follows: local submodules of projective modules are projective,
colocal factor modules of injective modules are injective. In particular: the
socle of any projective modules has support at some of the sinks, the top of
any injective module has support at some of the sources.
The Loewy length of a nicely tiered algebra with $n+1$ tiers is precisely $n+1$.
Of special interest is the following: 
{\it For a nicely  tiered algebra, any indecomposable module 
which is projective or injective is solid}
(an indecomposable module over an artin algebra is said to be {\it solid} provided its
socle series coincides with its radical series). 

 	\medskip
{\bf (8.1) Proposition.} {\it Let $\Lambda$ be nicely  tiered algebra with $n+1$ tiers. 
Assume that the following conditions are
satisfied for all indecomposable projective modules  $P, P'$ of Loewy length at least $3$:}
\item{\rm (P1)} {\it The module ${}_2P$ is a brick} (this means that any non-zero endomorphism
is an automorphism).
\item{\rm (P2)} {\it If $\Hom({}_2P,{}_2P') \neq 0$, then $P$ can be embedded into $P'.$}
	\smallskip
{\it Let $M$ be the direct sum of the modules ${}_tP$ with $P$ indecomposable projective
and $t\ge 2$ as well as the modules ${}_tQ$ with $Q$ indecomposable injective, 
and $t\ge 1$. Then $\End(M)$ has global dimension at most $n+2.$}
	\medskip

Remark. Note that the module $M$ considered in 8.1 is both a generator and
a cogenerator, thus the number $n+2$ is an upper bound for the representation
dimension of $\Lambda$. Here we encounter again a class of algebras $\Lambda$
where one now knows that $\repdim\Lambda \le \LL(\Lambda)+1.$

	\medskip
The proof of Proposition 8.1 follows the strategy of Iyama's proof [I] of the
finiteness of the representation dimension, as well as subsequent considerations
by Oppermann [O1], Corollary A1 (but see already [A] as well as
several joint papers with Dlab). It
relies on the following lemma shown in [R4]:
	\medskip
{\bf (8.2) Lemma.} {\it Let 
$$
 \emptyset = \Cal M_{-1} \subseteq \Cal M_0 \subseteq \cdots \subseteq \Cal M_{n+1} =
\Cal M
$$
be finite sets of indecomposable $\Lambda$-modules.
Let $M$ be a $\Lambda$-module with $\add M = \add \Cal M$ and $\Gamma = \End(M)$.

Assume that for any $N \in \Cal M_i$ 
there is a monomorphism 
$ u\:\alpha N \to N$
with $\alpha N\in \add \Cal M_{i-1}$ such that any radical map 
$\phi\:N' \to N$ with $N' \in \Cal M_i$
factors through $u$.

Then the global dimension of $\Gamma$ is at most $n+2$.}
	\medskip
Let us outline the proof of 8.2. We consider the indecomposable projective
$\Gamma$-modules $\Hom(M,N)$, where $N$ is indecomposable in $\Cal M$.
Assume that $N$ belongs to $\Cal M_i$ and not to
$\Cal M_{i-1}$, for some $i$. Since $N$ is not in $\Cal M_{i-1}$,
we see that $u$ is a proper monomorphism. Let 
$\Delta(N)$ be the cokernel of $\Hom(M,u)$, thus we deal with the exact
sequence
$$
 0 \to \Hom(M,\alpha N) @>\Hom(M,u)>> \Hom(M,N) \to \Delta(N) \to 0.
$$
We see that $\Delta(N)$ is the factor space of $\Hom(M,N)$
modulo those maps $M \to N$ which factor through $u$, thus through $\add \Cal M_{i-1}$.

The assumption that any radical map $\phi\:N' \to N$ with $N' \in \Cal M_i$
factors through $u$ means the following: if we consider $\Delta(N)$ as a
$\Gamma$-module, then it has one composition factor of the form 
$S(N) = \top \Hom(M,N)$,
all the other composition factors are of the form $S(N') = \top \Hom(M,N')$
with $N'$ an indecomposable module in $\Cal M$ which does not belong to $\Cal M_i$.

Since $\alpha N$ belongs to $\add \Cal M_{i-1}$, the projective
$\Gamma$-module $\Hom(M,\alpha N)$ is a direct sum of modules $\Hom(M,N')$
with $N'$ in $\Cal M_{i-1}.$ 

This shows that $\Gamma$ is left strongly quasi-hereditary with $n+2$ layers, thus
has global dimension at most $n+2$, according to [R4].
	\medskip

{\bf (8.3) Proof of 8.1.} Let $\Lambda$ be a nicely  tiered algebra with $n+1$ tiers
such that the conditions (P1) and (P2) are satisfied.
	\medskip
The conditions (P1) and  (P2) imply that corresponding properties
are satisfied for ${}_tP$ with $t > 2.$

\item{\rm (P1$_t$)} The module ${}_tP$ is a brick. 
\item{\rm (P2$_t$)} If $\Hom({}_tP,{}_tP') \neq 0$,
 then $P$ can be embedded into $P'.$

	\medskip
Proof: Let $f\:{}_tP \to {}_tP'$ be a non-zero homomorphism. Then also  $f|{}_2P$ is
non-zero, since otherwise $f$ would vanish at the tier $0$, but the socle of $P$
lives at the tier $0$. Thus, any non-zero endomorphism of ${}_tP$ yields
a non-zero endomorphism of ${}_2P$, by (P1) this is an isomorphism; but if $f|{}_2P$
has zero kernel, the same is true for $f\:{}_2P \to {}_2P$; thus $f$ is a mono endomorphism,
therefore an automorphism. This shows (P1$_t$). 
Similarly, if $f\:{}_tP \to {}_tP'$ is non-zero, then
also $f|{}_2P$ is non-zero, therefore 
$P$ can be embedded into $P'$ by (P2). 
	\medskip

We define sets of indecomposable modules $\Cal P^i, \Cal Q_i$
as follows:

Let $\Cal P^i$ be the set of modules ${}_tP$, 
where $P$ is indecomposable projective, $t\ge 2$
and $\LL(P)-t = i$. 
The modules in $\Cal P^i$ are indecomposable, according to condition (P1), see 8.1.
The non-empty sets $\Cal P^i$ are $\Cal P^0, \Cal P^1, \dots, \Cal P^{n-1}$; 
the modules in $\Cal P^0$
are the indecomposable projective modules which are not simple, those in $\Cal P^{n-1}$
are the modules of the form ${}_2P$ with $P$ generated at tier $n$. 

Let us collect some properties of the modules $N$ in $\Cal P^i.$ Such a module 
is generated at tier $g$ with $1 \le g  \le n-i$.
(Namely, if $P$ is indecomposable projective, then ${}_tP$ is of Loewy length $t$, thus
generated at $g = t-1$. Since $t \ge 2$, we have $g \ge 1$. 
Since $P$ is of length $l \le n+1$,
we have $t = l-i \le n+1-i$, thus $g = t-1 \le n-i.$)
The module $N$ lives in $[0,n-i]$ (since it is generated at tier 
$g \le n-i$), its socle lives at the vertices with tier $0$, and the 
Loewy length of such a module $N$ satisfies $2 \le \LL(N) \le n-i+1.$

Let $\Cal Q_i$ be the set of non-zero modules ${}_iQ$ with  $Q$ indecomposable injective
and $i \le \LL(Q).$ If $Q$ is cogenerated at tier $j$, where $0 \le j \le n$, then $\LL(Q)
= n-j+1$, thus $1 \le i \le n-j+1$ implies that $0 \le j \le  n-i+1$.
The non-empty sets $\Cal Q_i$ are $\Cal Q_1, \dots, \Cal Q_{n+1}.$
The modules in $\Cal Q_1$ are just all the simple modules. 

Since $i \ge 1$, the module ${}_iQ$ is a non-zero submodule of $Q$, thus has simple
socle and therefore is indecomposable. And again, 
we mention some additional properties for a module $N$ in $\Cal Q_i$. It 
is generated at tier $g$ with $i-1 \le g \le n$, it lives in $[i-1,n]$
and its Loewy length is precisely $i$.
	\medskip
Claim: {\it If $N$ is in $\Cal P^i$, then either $\rad N$ belongs to $\Cal P^{i+1}$,
or else $\rad N$ is semisimple and thus belongs to $\add\Cal Q_1.$
If $N$ in $\Cal Q_i$ with $i \ge 2$, then $\rad N$ belongs
to $\Cal Q_{i-1}$.}

Proof: Let $P$ be indecomposable projective of length $l$. If $t = l-i \ge 3,$
then $N = {}_tP$ belongs to $\Cal P^i$ and $\rad N = {}_{t-1}P$ belongs to $\Cal P^{i+1}.$
If $t = l-i = 2$, then $N = {}_tP$ has Loewy length $2$, thus $\rad N = {}_1P$ is semisimple,
and thus belongs to $\add \Cal Q_1.$

On the other hand, for $N = {}_iQ$, we have $\rad N = {}_{i-1}Q,$ and 
$n-i+1 < n-(i-1)-1$.
	\medskip

Now, let $\Cal M_i$ be the union of all the sets $\Cal P^j$ with $j \ge n+2-i$
as well as the sets $\Cal Q_j$ with $j\le i.$
Thus,
$$
\align
  \Cal M_{n+2} &= \quad \Cal P^{0}\quad \qquad\quad\  \ \cup \Cal M_{n+1} \cr
  \Cal M_{n+1} &= \quad  \Cal P^{1}\quad\cup \Cal Q_{n+1}  \cup \Cal M_n\cr
  \Cal M_n     &= \quad \Cal P^{2}\quad\cup \Cal Q_n  \quad \cup \Cal M_{n-1}\cr
    \cdots &  \cr
  \Cal M_i &= \Cal P^{n+2-i}\cup \Cal Q_i \quad \cup \Cal M_{i-1} \cr
    \cdots &  \cr
  \Cal M_3 &= \ \Cal P^{n-1}\ \cup \Cal Q_3 \quad \cup \Cal M_2 \cr
  \Cal M_2 &= \qquad\qquad\  \Cal Q_2  \quad \cup \Cal M_1  \cr
  \Cal M_1 &= \qquad\qquad\  \Cal Q_1  \cr
  \Cal M_0 &= \quad \ \emptyset 
\endalign
$$
As we have shown: If $N$ belongs to $\Cal M_i$ for some $i$, then 
$\rad N$ is in $\add \Cal M_{i-1}.$ Thus, let $\alpha N = \rad N$ and $u\:\alpha N \to N$
the inclusion map. We want to verify that we can apply Lemma 8.2.

Thus, we have to show that for $N \in \Cal M_i$ any non-zero radical map 
$f\:N' \to N$ with $N' \in \Cal M_i$ maps into $\rad N$. 
We can assume that $N \notin \Cal M_{i-1},$ thus $N$ belongs either to $\Cal P^{n+2-i}$
or to $\Cal Q_i$.

First, let us assume that $N \in \Cal Q_i$, thus $N$ has Loewy length $i$.
If $N'\in \Cal P^{n+2-j}$ with $j \le i$, then,  $N'$ has Loewy length 
at most $n-(n+2-j)+1 = j-1 \le i-1$. 
Similarly, if $N' \in \Cal Q_j$, with $j < i$, then the Loewy length of $N'$
is at most $i-1$. In both cases, we see that the image of any map $f\:N' \to N$
lies in $\rad N$. Thus, it remains to consider the case that $N' \in \Cal Q_i$,
so that $N'$ has also Loewy length $i$. Now $N'$ has a simple socle. If $f$
vanishes on the socle, then again the image of $f$ has socle length at most $i-1$
and thus lies in $\rad N$. 
If $f$ does not vanish on the socle, then $f$ is a monomorphism.
But $N'$ is relative injective in the subcategory of all modules of Loewy length at
most $i$, thus $f$ is a split monomorphism, thus not a radical morphism.
	\medskip

Second, we assume that $N \in \Cal P^{n+2-i}.$
First, consider the case that $N'\in \Cal Q_j$ with $j \le i$.
Now the socle of $N$ lives at tier $0$, thus the image of $f$ (and therefore
$N'$ itself) must have a composition factor at tier $0$. This shows that
$N' = {}_jQ$ with $Q$ the injective envelope of a simple at tier $0$ and
that $f$ is injective. Assume that the image of $f$ does not lie in $\rad N$,
then the Loewy length of $N$ has to be equal to $j$. But ${}_jQ$ is relative injective
in the subcategory of all modules of Loewy length at
most $j$, thus $f\:N' \to N$ is a split mono, thus not a radical map. 

Finally, there is the case that $N'\in \Cal P^{n+2-j}$ with $j \le i.$ 
Let $N = {}_tP$ and $N' = {}_{t'}P'$ with $P$ of Loewy length $l$ and $P'$ of
Loewy length $l'$. 
If $t' < t$, then $f\:{}_{t'}P' \to {}_tP$ maps into the radical of ${}_tP.$
If $t' > t$, then $\Hom({}_{t'}P',{}_tP) = 0$, since ${}_{t'}P'$ is generated at the
tier $t'$, and ${}_tP$ lives at the tiers $[0,t].$ Thus, we can assume that $t' = t.$ 
Since $j \le i$, we see that $\LL(P') = n+2-j+t \ge n+2-i+t = \LL(P).$
If $\Hom({}_tP',{}_tP) \neq 0$, then $P'$ can be embedded
into $P$, according to condition (P1), thus $\LL(P') \le \LL(P)$
and therefore $\LL(P') = \LL(P).$ But if $P'$ is isomorphic to a submodule of $P$ and both
have the same Loewy length, then  $P'$ and $P$ are isomorphic 
and therefore also ${}_tP'$ and ${}_tP$ are isomorphic. 
But then we  use (P1$_t$) in order to see that any non-zero homomorphism ${}_tP' \to {}_tP$
is an isomorphism. This contradicts the assumption that there is a
non-zero radical map ${}_tP' \to {}_tP$.
	\medskip
Remark.  One should be aware that the classes $\Cal P^j$ and $\Cal Q_i$
are not necessarily disjoint. A typical example is the fully commutative square
$$
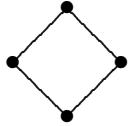
{\beginpicture
\setcoordinatesystem units <.6cm,.6cm>
\plot 0 1  1 0  2 1  1 2  0 1  /
\multiput{$\bullet$} at 0 1  1 0  2 1  1 2 /
\endpicture}
$$
(say with arrows pointing downwards). This is a nicely  tiered algebra with 3 tiers.
There is an indecomposable module $P$ which is projective-injective, it belongs both to
$\Cal P^0$ and to $\Cal Q_3$.

	\medskip
{\bf (8.4)} Let us add some examples 
of nicely  tiered algebras which do not satisfy the conditions (P1), (P2), respectively.
Again, we present the quivers by just indicating the corresponding
edges; all the arrows are supposed to point downwards. 
$$
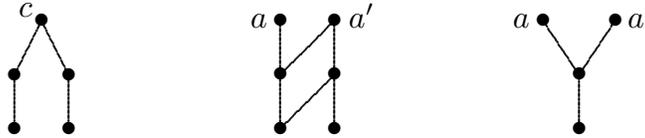
{\beginpicture
\setcoordinatesystem units <.6cm,.6cm>
\put{\beginpicture

\plot 0 0  0 1  0.5 2  1 1  1 0 /
\multiput{$\bullet$} at 0 0  0 1  0.5 2  1 1  1 0  /
\put{$c$} at 0.2 2.2
\endpicture} at 0 0
\put{\beginpicture
\setcoordinatesystem units <.6cm,.6cm>
\plot 1 1  0 0  0 2  /
\plot 0 1  1 2  1 0 /
\plot 1 2  1 2 /
\multiput{$\bullet$} at 0 0  0 1  0 2  1 0  1 1  1 2 /
\put{$a'$} at 1.5 2.1
\put{$a$} at -.4 2
\endpicture} at 5 0
\put{\beginpicture
\setcoordinatesystem units <.4cm,.6cm>
\plot 0 2  1 1  2 2 /
\plot 1 1  1 0 /
\multiput{$\bullet$} at 1 0  1 1  0 2  2 2 /
\put{$a'$} at 2.7 2.1
\put{$a$} at -.6 2
\endpicture} at 10 0
\endpicture}
$$
In the example left, ${}_2P(c)$ is decomposable. In the middle example, 
we consider the
path algebra of the quiver with the commutativity relation. Then both
${}_2P(a)$ and 
${}_2P(a')$ are indecomposable. We see that 
$\Hom({}_2P(a),{}_2P(a')) \neq 0,$ but $P(a)$ cannot be embedded into $P(a').$
On the right, we see a further example where the condition (P1), but not
the condition (P2) is satisfied.

	\bigskip\medskip
{\bf 9. Tensor products of path algebras of bipartite quivers.}
	\medskip
Recall that a finite quiver is said to be {\it bipartite} if and only if every
vertex is a sink or a source. Thus, a quiver $Q$ is bipartite if and only if
its path algebra is a finite dimensional algebra with radical square zero.
	\medskip
{\bf (9.1) Theorem.} {\it 
Let $\Lambda_1,\dots\Lambda_n$ be path algebras of bipartite quivers. 
Then the algebra
$\Lambda = \Lambda_1\otimes_k\cdots\otimes_k\Lambda_n$ 
has representation dimension at most $n+2$.}
	\medskip
For the proof, we want to use Proposition 8.1. Of course, we can assume that all the
algebras $\Lambda_i$ are connected and not simple, thus tiered with precisely $2$
tiers. In order to show that
$\Lambda$ is tiered with $n+1$ tiers, we use induction and the following general result:
	\medskip
{\bf (9.2)} {\it The tensor product of nicely  tiered algebras with $n_1+1$ and $n_2+1$ tiers
respectively is nicely  tiered with $n_1+n_2+1$ tiers.}
	\medskip
Proof. Let $\Lambda_1$ and $\Lambda_2$ be nicely  tiered algebras with $n_1+1$ and $n_2+1$
tiers, respectively. Let $Q^{(1)}$ be the quiver of $\Lambda_1$ and 
$Q^{(2)}$ that of $\Lambda_2$.
Then the quiver of $\Lambda_1\otimes_k \Lambda_2$ is $Q = Q^{(1)}\otimes Q^{(2)}$,
this is the quiver with vertex set $Q^{(1)}_0\times Q^{(2)}_0$, and with arrow set
$(Q^{(1)}_1\times Q^{(2)}_0) \cup (Q^{(1)}_0\times Q^{(2)}_1)$; 
here, given an arrow $\alpha_1\:x_1\to y_1$ in $Q^{(1)}$,
and a vertex $z_2$ of $Q^{(2)}$, there is the arrow $\alpha_1z_2\:x_1z_2
\to y_1z_2$,
and similarly, given a vertex $x_1$ of $Q^{(1)}$ and an arrow $\beta_2\:y_2\to z_2$
in $Q^{(2)}$, there is the arrow $x_1\beta_2\:x_1y_2 \to x_1z_2.$
(When writing down an element  of a product $U\times V$, we just write 
$uv$ instead of $(u,v)$, for $u\in U$ and $v\in V$.)

For example, in case we consider the tensor product of two copies of the
Kronecker algebra, say with quivers $Q^{(1)}$ and $Q^{(2)}$, 
we obtain the following quiver $Q = Q^{(1)}\otimes Q^{(2)}$:
$$
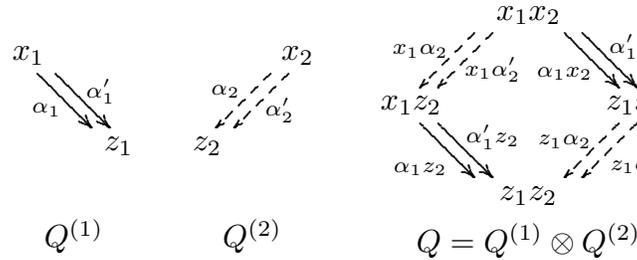
{\beginpicture
\setcoordinatesystem units <1cm,1cm>
\put{\beginpicture
\put{$x_1$} at 1 2
\put{$z_1$} at 2 1
\arr{1.3 1.8}{1.9 1.2}
\arr{1.1 1.8}{1.7 1.2}
\put{$\ssize \alpha'_1$} at 1.8 1.6
\put{$\ssize \alpha_1$} at 1.2 1.4
\put{$Q^{(1)}$} at 1.5 0
\put{} at 1.5 2.5
\endpicture} at -5 0
\put{\beginpicture
\put{$x_2$} at 2 2
\put{$z_2$} at 1 1
\setdashes <1mm>
\arr{1.9 1.8}{1.3 1.2}
\arr{1.7 1.8}{1.1 1.2}
\put{$\ssize \alpha_2$} at 1.2 1.6
\put{$\ssize \alpha'_2$} at 1.8 1.4
\put{$Q^{(2)}$} at 1.5 0
\put{} at 1.5 2.5

\endpicture} at -3 0

\put{\beginpicture
\put{$x_1x_2$} at 1 2
\put{$x_1z_2$} at -.3 1
\put{$z_1x_2$} at 2.2 1
\put{$z_1z_2$} at 1 0

\arr{1.6 1.8}{2.2 1.2}
\arr{0 0.8}{0.6 0.2}
\arr{1.4 1.8}{2 1.2}
\arr{-.2 0.8}{0.4 0.2}
\setdashes <1mm>
\setdashes <1mm>
\arr{0.6 1.8}{0 1.2}
\arr{2.2 0.8}{1.6 0.2}
\arr{0.4 1.8}{-.2 1.2}
\arr{2 0.8}{1.4 0.2}

\put{$\ssize x_1 \alpha_2$} at -.2 1.6
\put{$\ssize x_1 \alpha'_2$} at 0.6 1.4

\put{$\ssize \alpha_1x_2$} at 1.4 1.35
\put{$\ssize \alpha'_1 x_2$} at 2.2 1.65

\put{$\ssize \alpha_1 z_2$} at -.2 0.3
\put{$\ssize \alpha'_1 z_2$} at .6 0.65
\put{$\ssize z_1 \alpha_2$} at 1.4 0.6
\put{$\ssize z_1 \alpha'_2$} at 2.2 0.35

\put{$Q = Q^{(1)}\otimes Q^{(2)}$} at 1 -.5

\endpicture} at 0 0

\endpicture}
$$
Here, the arrows of $Q^{(1)}$ as well as those of $Q$ which belong to
$Q^{(1)}_1\times Q^{(2)}_0$ are shown as solid arrows, those of 
$Q^{(2)}$ as well as those of $Q$ which belong to
$Q^{(1)}_0\times Q^{(2)}_1$ are shown as dashed ones.

Now suppose that $Q^{(1)},Q^{(2)}$ are nicely  tiered, with $n_1+1$ and $n_2+1$ tiers,
and tier functions $l_1,l_2$ respectively. Given vertices $x_1\in Q^{(1)}_0$
and $x_2\in Q^{(2)}_0$, define $l(x_1x_2) = l_1(x_1)+l_2(x_2)$. This defines a function
$Q^{(1)}_0\times Q^{(2)}_0 \to \Bbb Z$ with values in the interval $[0,n+n'].$
For a sink $x_1x_2$ of $Q^{(1)}\otimes Q^{(2)}$, we have $l(x_1x_2) = 0$, for a source 
$(x_1x_2)$ of $Q^{(1)}\otimes Q^{(2)}$, we have
$l(x_1x_2) = n_1+n_2$ 
and given an arrow of $Q^{(1)}\otimes Q^{(2)}$, the value of $l$ decreases by $1$.
This shows that we obtain a tier function for a nicely tiered quiver. 

Finally, we have to note that the relations of $\Lambda_1\otimes_k \Lambda_2$ are obtained from
the relations of $\Lambda_1$ and $\Lambda_2$ and adding commutativity relations;
thus, if $\Lambda_1$ and $\Lambda_2$ are defined by using only commutativity relations,
the same is true for $\Lambda_1\otimes_k \Lambda_2$.

	\medskip
{\bf (9.3)} {\it 
Let $\Lambda_1,\dots\Lambda_n$ be path algebras of bipartite quivers. 
Then $\Lambda = \Lambda_1\otimes_k\cdots\otimes_k\Lambda_n$ 
satisfies the conditions {\rm (P1)} and {\rm(P2)}.}
	\medskip
Proof. 
Let us introduce some notation concerning 
tensor products $\Lambda = \Lambda_1\otimes_k \cdots \otimes_k\Lambda_n$, where any 
$\Lambda_i$ is the path algebra of a finite directed quiver $Q^{(i)}.$
The quiver $Q$ of $\Lambda$ is given as follows: The set of vertices is the set
$Q^{(1)}_0\times\cdots\times Q^{(n)}_0$, an element of this set will be denoted 
by $x = x_1x_2\cdots x_n$ with $x_i\in Q^{(i)}_0$ for $1\le i \le n.$ 
Given such a vertex $x$, we are interested in the corresponding indecomposable
projective module $P(x).$ 
	\medskip
Let $W(x_i,y_i)$ be the set of paths in $Q^{(i)}$ starting in $x_i$ and ending in $y_i$,
this may be considered as a basis of $P(x_i)_{y_i}$ and therefore we may take as
basis of $P(x)_y$, where $x,y$ are vertices of $Q$,
the product set $W(x,y) = W(x_1,y_1)\times\cdots\times W(x_n,y_n)$; 
we call this the {\it path basis} of $P(x)$.
In particular, we see: a vertex $y$ belongs to the support of $P(x)$
if and only if there are paths starting at $x_i$ and ending in $y_i$,
for $1\le i \le n.$ 

Now assume that $Q^{(i)}$ a bipartite, thus any path in $Q^{(i)}$ is of length at most $1$, thus either a vertex or an arrow. We want to 
describe the representation ${}_2P(x)$ for any
vertex $x = x_1\dots x_n$ of $Q$. 
Note that the support quiver of ${}_2P(x)$ will again be bipartite. 
We can assume that $t$ of the vertices $x_i$
are sources, and the remaining ones sinks. Thus, up to a permutation we can assume
that $x = x_1\cdots x_n$ with sources $x_i$ for $1\le j\le t$ 
and sinks $x_j = z_j$ for $t+1\le j \le n.$ 
The support $S$ of the socle of $P(x)$ consists of the
vertices $z = z_1\cdots z_n$ where  $z_i$
is a sink in the quiver $Q^{(i)}$ such that there is a path from $x_i$ to $z_i$,
for any $1\le i\le t.$ 

Given an $n$-tuple $u_1u_2\cdots u_n$ where the $u_i$ are elements of some 
sets (say of vertices or arrows of some quivers), and $v_j$ is a further element, then
we denote by $u[ v_j = u_1\dots u_{j-i}v_ju_{j+1}\cdots u_n$ the element
obtained from $u$ by replacing its entry at the position $j$ by $v_j.$

Using this notation, the vertices in the support of ${}_2P(x)$ with
tier number $1$ are of the
form $y = z[x_j$ and the arrows of the support
are of the form $z[\alpha_j$, always with $z \in S$, and with arrows
$\alpha_j\:x_j \to z_j,$  and $1\le j \le t.$ 

If we are interested in the structure of $P(x),$ we may assume that all the vertices
$x_j$ are sources, thus that $t=n$ 
(namely, if for example $x_n$ is a sink, then $P(x_n)$ is one-dimensional
and thus $P(x) = P(x_1\cdots x_{n-1})\otimes_k P(x_n)$ can be identified with
the $(\Lambda_1\otimes_k\cdots\otimes_k\Lambda_{n-1})$-module $P(x_1\cdots x_{n-1})$).

Given a vertex $x$ of $Q$, we may look at the coefficient quiver $\Theta (x)$ 
of $P(x)$  with respect to its path basis (for the definition, see [R2]).
If we look at our example of the tensor product of two copies of the
Kronecker algebra, and consider the unique source $x = x_1x_2$ of
$Q$, then the coefficient quiver $\Theta (x)$ 
of $P(x)$ with respect to the path basis looks as shown on the right:
$$
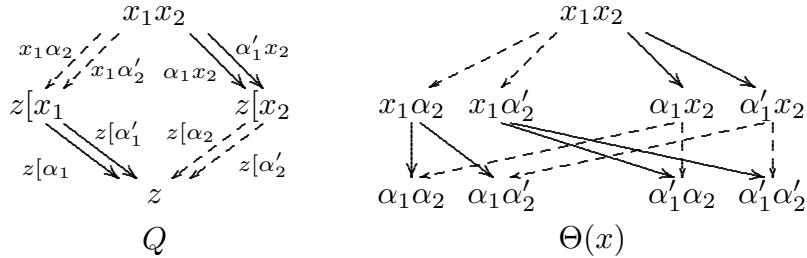
{\beginpicture
\setcoordinatesystem units <1cm,1cm>
\put{\beginpicture
\put{$x_1x_2$} at 1 2
\put{$z[x_1$} at -.3 1
\put{$z[x_2$} at 2.2 1
\put{$z$} at 1 0

\arr{1.6 1.8}{2.2 1.2}
\arr{1.4 1.8}{2 1.2}
\arr{0 0.8}{0.8 0.2}
\arr{-.2 0.8}{0.6 0.2}
\setdashes <1mm>
\setdashes <1mm>
\arr{0.6 1.8}{0 1.2}
\arr{0.4 1.8}{-.2 1.2}
\arr{2.2 0.8}{1.4 0.2}
\arr{2 0.8}{1.2 0.2}

\put{$\ssize x_1 \alpha_2$} at -.2 1.6
\put{$\ssize x_1 \alpha'_2$} at 0.6 1.4

\put{$\ssize \alpha_1x_2$} at 1.4 1.35
\put{$\ssize \alpha'_1 x_2$} at 2.2 1.65

\put{$\ssize z[\alpha_1$} at -.2 0.3
\put{$\ssize z[\alpha'_1$} at .6 0.7
\put{$\ssize z[\alpha_2$} at 1.4 0.7
\put{$\ssize z[\alpha'_2$} at 2.2 0.35

\put{$Q$} at 1 -.5

\endpicture} at 0 0

\put{\beginpicture
\put{$x_1x_2$} at 1 2

\put{$x_1\alpha_2\strut$} at -1 1
\put{$x_1\alpha'_2\strut$} at 0 1
\put{$\alpha_1x_2\strut$} at 2 1
\put{$\alpha'_1x_2\strut$} at 3 1

\put{$\alpha_1\alpha_2\strut$} at -1 0
\put{$\alpha_1\alpha'_2\strut$} at 0 0
\put{$\alpha'_1\alpha_2\strut$} at 2 0
\put{$\alpha'_1\alpha'_2\strut$} at 3 0

\arr{1.6 1.8}{2.8 1.2}
\arr{1.4 1.8}{2 1.2}
\arr{-1 0.8}{-1 0.2}
\arr{-.9 0.8}{-.1 0.2}
\arr{0 0.8}{1.9 0.2}
\arr{0.1 0.8}{2.9 0.2}

\setdashes <1mm>
\arr{0.6 1.8}{0 1.2}
\arr{0.4 1.8}{-.8 1.2}
\arr{1.9 0.8}{-.9 0.2}
\arr{2.9 0.8}{.1 0.2}
\arr{2 0.8}{2 0.2}
\arr{3 0.8}{3 0.2}



\put{$\Theta (x)$} at 1 -.5

\endpicture} at 4.9 0
\endpicture}
$$
On the left, we present again the quiver of $Q = Q^{(1)}\otimes Q^{(2)}$,
but now using the notation $z[?$ for the vertices with tier number $1$ as well
as the arrows ending in $z$.
	\medskip
Now we are going to look at ${}_2P(x)$
Let  $z\in S$ and take an arrow $\alpha_j\:x_j \to z_j$,
let $y = z[x_j$, this is a vertex with tier number $1$. 
The vector spaces $P(x)_y, P(x)_z$ and the linear map $P(x)_{z[\alpha_j}$
are given as follows: 
Since $z$ belongs to $S$, there is an arrow $x_i\to z_i$ for any $i$
and $W(x,z)$ is a basis of $P(x)_z$ (note that here $W(x_i,z_i)$ is the set of
arrows $x_i \to z_i$ for all $i$). 
Similarly, for $y = z[x_j$, the space $P(x)_y$ has as a basis the set of
elements of the form $\alpha[x_j$ with $\alpha\in W(x,z)$, and the linear map
$z[\alpha_j\: P(x)_y \to P(y)_z$ sends $\alpha[x_j$ to $ \alpha[\alpha_j.$

In the coefficient quiver ${}_2\Theta(x)$ of ${}_2P(x)$, 
any sink $\alpha$ is the end point of precisely
$n$ arrows, namely the arrows labeled $z[\alpha_j\:z[x_j \to z$,
where $\alpha_j\in W(x_j,z_j).$ 
It follows that the one-dimensional vector space $k\alpha$ generated
by $\alpha$ is the intersection of the images of the maps 
$$
 k\alpha = \bigcap_{j=1}^n \Im\left(
 z[\alpha_j\: P(x)_{z[x_j} \longrightarrow P(x)_z \right).
$$
As a consequence, any endomorphism of ${}_2P(x)$ will map the element
$\alpha$ of $P(x)_z$ onto 
a multiple of $\alpha$.
	
We claim that {\it the coefficient quiver ${}_2\Theta(x) $ 
of ${}_2P(x)$ with respect
to the path basis is connected}.
Namely, given a sink $\alpha = \alpha_1\cdots \alpha_n$ of ${}_2\Theta(x) $,
and any arrow $\alpha'_j\:x_j \to z'_j$ in $Q^{(j)}$ different from $\alpha_j$, 
there is a path of length $2$ in ${}_2\Theta(x) $ starting in $\alpha$ and ending in
$\alpha[\alpha'_j$, namely
$$
{\beginpicture
\setcoordinatesystem units <1cm,1cm>
\put{$\alpha$} at 0 0
\put{$\alpha[x_j$} at 1 1
\put{$\alpha[\alpha'_j$} at 2 0
\arr{0.8 0.8}{0.2 0.2}
\arr{1.2 0.8}{1.8 0.2}
\endpicture}
$$
Thus, given two sinks of ${}_2\Theta(x) $, say $\alpha = \alpha_1\cdots \alpha_n$
and $\alpha' = \alpha'_1\cdots \alpha'_n$, we may replace successively
$\alpha_j$ by $\alpha'_j$ and obtain a path of length at most $2n$ starting
in $\alpha$ and ending in $\alpha'.$
	\medskip
In order to deal with the conditions (P1) and (P2), 
we consider maps $f\:{}_2P(x) \to M$ with $M = P(x')$ for some
vertex $x'$. Actually, the essential property of $M$ which we will need
is that all the maps used in $M$ are injective. Thus, let 
$\Cal M$ be the set of $\Lambda$-modules $M$ such that all the maps
used are injective. Clearly, all the indecomposable projective
$\Lambda$-modules, and even all their submodules belong to $\Cal M$.

Consider $f\:{}_2P(x) \to M$ with $M \in \Cal M$.
We show: {\it Given any arrow $\alpha[x_j \to \alpha$ in the coefficient quiver,
then $f(\alpha) = 0$ if and only if $f(\alpha[x_j) = 0.$}

Proof: Clearly, if $f(\alpha[x_j) = 0$, then also $f(\alpha) = 0$, since
$\alpha$ is a multiple of $\alpha[x_j$. Thus, conversely, let us assume that
$f(\alpha) = 0$. There is the following commutative diagram
$$
\CD
 P(x)_{z[x_j} @>f_{z[x_j}>> M_{z[x_j} \cr
  @VP(x)_{z[\alpha_j}VV        @VVM_{z[\alpha_j}V   \cr
 P(x)_z @>f_z>> M_z,
\endCD  
$$
with $\alpha[x_j$ being sent by the left vertical map to $\alpha$.
Since the right vertical map is injective, the vanishing of $f_z(\alpha)$
implies that also $f(\alpha[x_j) = 0.$

As a consequence of the connectivity of the 
coefficient quiver of ${}_2P(x)$ we conclude: {\it If $f(\alpha) = 0$ for some $\alpha$,
then $f = 0.$}

	\smallskip
Condition (P1). Let $f$ be an endomorphism of ${}_2P(x)$. Assume that $f$
is not a monomorphism. Since $f$ maps any basis vector $\alpha$ of the
socle of ${}_2P(x)$ onto a multiple of itself, we see that there has to be 
such a basis element $\alpha$ with $f(\alpha) = 0.$ But then $f = 0.$
 	\smallskip
Condition (P2). 
Assume there is given a non-zero homomorphism $f\:{}_2P(x) \to M$
with $M \in \Cal M$. Assume that $f(P(x)_{z[x_j}) = 0.$ Of course, then
also $f(P(x)_z) = 0.$ Thus, it follows again that $f = 0$.

Thus we see: if there is a non-zero homomorphism $f\:{}_2P(x) \to P(x')$, then
all the elements $z[x_j$ with $1\le j \le t$ 
are in the support of $P(x')$, and therefore $x'_j = x_j.$
It follows that $P(x)$ is a submodule of $P(x').$
	\medskip
Remark: An alternative way for proving Theorem 9.1 is as follows. First,
consider the special case where none of the quivers $Q^{(i)}$ has multiple arrows.
Under this assumptions all the indecomposable projective $\Lambda$-modules are thin,
thus we do not have to worry about bases. 
The general case can then be obtained from this special case using covering theory. 
	\medskip
{\bf (9.4) Corollary.} {\it 
Let $\Lambda_1,\dots\Lambda_n$ be path algebras of representation-infinite 
bipartite quivers. The algebra
$\Lambda = \Lambda_1\otimes_k\cdots\otimes_k\Lambda_n$ 
has representation dimension precisely $n+2$.}
	\medskip
Proof. This follows directly from the inequalities 7.5 and 
and 9.1.
	\bigskip
The special case of the $2$-fold tensor power of the Kronecker algebra 
has been exhibited by Oppermann in [O2],
when he considered one-point extensions of wild algebras. This example was
the starting point for our investigation.

	\bigskip\medskip
{\bf References}

{\frenchspacing
	\medskip

\item{[APT]} Assem, I., Platzek, M.I., Trepode, S.: 
  On the representation dimension
  of tilted and laura algebras, J. Algebra 296 (2006), 426-439.

\item{[A]} Auslander, M.: The representation dimension of artin algebras. Queen Mary
   College Mathematics Notes (1971). Republished in {\it Selected works of Maurice
Auslander.} Amer. Math. Soc., Providence 1999.

\item{[AB]} Auslander, M., Bridger, M.: Stable Module Theory. Memoirs Amer.
   Math. Soc. 94 (1969). 

\item{[AR1]} M. Auslander, I. Reiten: Stable equivalence of dualizing $R$-varieties:
    III. Dualizing $R$-varieties stably equivalent to hereditary dualizing
    $R$-varieties.  Adv. in Math. 17 (1975), 122-142.

\item{[AR2]} M. Auslander, I. Reiten: Stable equivalence of dualizing $R$-varieties:
    V. Artin algebras stably equivalent to hereditary algebras. Adv. in Math. 17
    (1975), 167-195.

\item{[Bs]} Bass, H.: Finitistic dimension and homological generalization of
semiprimary rings. Trans Amer. Math. Soc. 95 (1960), 466-488.

\item{[Bt]} Bautista, R.: On algebras of strongly unbounded representation
type. Comment. Math. Helv. 60 (1985), 392-399.

\item{[Bi]} BIREP-Workshop 2008: The Representation Dimension of Artin Algebras.
  \newline http://www.math.uni-bielefeld.de/$\sim$sek/2008/

\item {[Bo]} Bongartz, K.: Indecomposables are standard. Comment.
Math. Helv. 60 (1985), 400-410.

\item {[Bo2]} Bongartz, K.: 
   Indecomposables live in all smaller lengths. Preprint. arXiv:0904.4609

\item{[CE]} Cartan, H., Eilenberg, S.: Homological Algebra. Princeton University Press,
  (1956).

\item{[CP]} F. Coelho, M. I. Platzeck: On the representation dimension of some classes
  of algebras. J. Algebra 275 (2004), 615-628.

\item{[EHIS]} K. Erdmann, Th. Holm, O. Iyama, J. Schr\"oer: Radical embeddings 
  and representation dimension. Adv. in Math. 185 (2004), 159-177.

\item{[G]} Gabriel, P.: The universal cover of a representation-finite algebra.
 In: Representations of Algebras. Springer LNM 903 (1981), 68-105.

\item{[HU]} Happel, D., Unger, L.: Representation dimension and tilting.
Journal Pure Appl. Algebra 215 (2011), 2315-2321.

\item{[IT]} Igusa, K., Todorov, G.: On the finitistic global dimension conjecture for Artin   algebras. Fields Institute Communications, vol. 45, 2005, 201-204.  

\item{[I]}
 Iyama, O.: Finiteness of representation dimension. Proc\. Amer\. Math\. Soc.
 131 (2003), 1011-1014.

\item{[J]} Jans, J.P.: Rings and Homology. Holt, Rinehart and Winston. New York 1964.

\item{[L]} Leszczy\'nski, Z.: On the representation type of tensor product algebras,
  Fund. Math. 144 (1994), 143-161.

\item{[LS]} Leszczy\'nski, Z., Skowro\'nski, A.: Tame tensor products of algebras,
   Colloq. Math. 98 (2003), 125-145.

\item{[O1]} Oppermann, S.: 
 Lower bounds for Auslander's representation dimension
 Duke Math. J. 148 (2009), 211-249.

\item{[O2]} Oppermann, S.: 
 Wild algebras have one-point extensions of representation dimension at least four
 J. Pure Appl. Algebra 213 (2009), no. 10, 1945-1960. 

\item{[O3]} Oppermann, S.; Representation dimension of quasi-tilted algebras
   J. Lond. Math. Soc. (2) 81 (2010), 435-456.

\item{[O4]} Oppermann, S.: Private communication, May 6, 2008.

\item{[R1]} Ringel, C. M.: Tame algebras and integral quadratic forms. 
 Springer LNM 1099.

\item{[R2]} Ringel, C.M.: Exceptional modules are tree modules. 
  Lin. Alg. Appl. 275-276 (1998) 471-493

\item{[R3]} Ringel, C.M.: The torsionless modules of an artin algebra. Selected Topics
 WS 2007/8. 
 http://www.math.uni-bielefeld.de/~ringel/opus/torsionless.pdf 

\item{[R4]} Ringel, C.M.:
 Iyama's finiteness theorem via strongly quasi-hereditary algebras.
 Journal of Pure and Applied Algebra 214 (2010) 1687-1692.

\item{[R5]} Ringel, C.M.: Indecomposables live in all smaller lengths. 
    Bulletin of the London Mathematical Society 2011; 
    doi: 10.1112/blms/bdq128 

\item{[R6]} Ringel, C.M.: The minimal representation-infinite algebras which are special     biserial. To appear in: Representations of Algebras and Related Topics, 
   EMS Series of Congress Reports, European Math. Soc. Publ. House, Z\"urich, 2011

\item{[Rq]} Rouquier, R.:
 Representation dimension of exterior algebras, Inventiones Math. 165 (2006), 357-367. 

\item{[X1]} Xi, Ch.: On the representation dimension of finite dimensional algebras.
 Journal of Algebra, 226, (2000), 332-346.

\item{[X2]} Xi, Ch.: Representation dimension and quasi-hereditary algebras.
  Advances in Math\. 168 (2002), 193-212.
\par}

	\bigskip\medskip
{\parindent=0truecm
\parskip=-2pt

\rmk Fakult\"at f\"ur Mathematik, Universit\"at Bielefeld \par
POBox 100\,131,\par D-33\,501 Bielefeld, Germany \par
e-mail: \ttk ringel\@math.uni-bielefeld.de \par}

\bye